\documentclass[11pt]{amsart}
\usepackage{graphicx}
\usepackage{amssymb}
\usepackage{amsmath}
\usepackage{amsthm}
\usepackage{epstopdf}
\usepackage[mathscr]{eucal}
\usepackage{verbatim}
\usepackage{amssymb}
\usepackage{epsfig}
\usepackage{enumerate}

\newcommand{\N}{\mathbb{N}}

\newcommand{\Z}{\mathbb{Z}}

\newtheorem{theorem}{{\bf Theorem}}
\newtheorem{lemma}{\noindent {\bf Lemma}}
\newtheorem{proposition}{\noindent {\bf Proposition}}
\newtheorem{corollary}{\noindent {\bf Corollary}}

\newtheorem{claim}{\noindent {\bf Claim}}
\newtheorem{fact}{\noindent {\bf Fact}}

\newtheorem{problem}{\noindent{\bf Problem}}

\def\endproof{\hfill {\kern 6pt\penalty 500
\raise -0pt\hbox{\vrule \vbox to5pt {\hrule width 5pt
\vfill\hrule}\vrule}}}

\begin{document}

\title[The morphology of infinite tournaments]{The morphology of infinite tournaments. Application to  the growth of their profile}
\author[Y.Boudabbous]{Youssef Boudabbous}
\address{ D\'epartement de Math\'ematiques, Facult\'e des Sciences,  Universit\'e de Sfax,   B.P.802, 3018 Sfax, Tunisie}
\email{youssef.boudabbous@yahoo.fr}
\author [M.Pouzet]{Maurice Pouzet } \thanks{This work was done under the auspices of the Intas programme Universal algebra theory and   supported by CMCU Franco-Tunisian "Outils math\'ematiques pour l'informatique"}
\address{UFR de Math\'ematiques, Universit\'e Claude-Bernard  Lyon 1, 
43, Bd. du 11 Novembre 1918,
69622 Villeurbanne, France}
\email {pouzet@univ-lyon1.fr}


\dedicatory{Dedicated to Michel  Deza  at the occasion of his
$65^{th}$ birthday.}

\date{\today}

\begin{abstract} A tournament is \emph{acyclically indecomposable} if no acyclic autonomous set of vertices has more than one element. We identify twelve infinite acyclically indecomposable tournaments and prove that every infinite acyclically indecomposable tournament contains a subtournament isomorphic to one of these  tournaments. The {\it profile} of a tournament $T$ is the function $\varphi_T$ which counts for each integer $n$ the number $\varphi_T(n)$ of  tournaments induced by $T$ on the  $n$-element subsets of $T$, isomorphic tournaments being identified. As a corollary of the result above we deduce that the growth of $\varphi_T$ is either polynomial, in which case  $\varphi_T(n)\simeq an^k$, for some positive real $a$, some non-negative integer $k$, or  as fast as some exponential.   \end{abstract}
\maketitle
\section{Introduction and presentation of the results}
An important chapter of the theory of graphs is about the decompositions  of graphs into simpler subgraphs. A wealth of results has been obtained  along the lines pioneered by   T. Gallai \cite {gallai},  \cite{maffray} with the notion of indecomposable graph (see \cite{ehrenfeucht} for an example). This paper is about tournaments. We study acyclically indecomposable tournaments, objects  introduced by Culus and Jouve in 2005 \cite{culus-jouve1}.  A consequence of our study is the existence of a gap in the growth rate of  the profile of tournaments.   

\subsection{Lexicographical sums of acyclic tournaments and acyclically indecomposable tournaments}
A  tournament is \emph{acyclically indecomposable} if no acyclic autonomous set of vertices has more than one element. 
Our first result is quite elementary:
\begin{theorem} \label{acyclicallyindecomposable1} Every tournament  $T$ decomposes into a lexicographical sum of acyclic tournaments indexed by an acyclically indecomposable tournament. The decomposition is unique and up to an isomorphism, this acyclically indecomposable tournament is unique. 
\end{theorem}
The blocks of the decomposition are the \emph{acyclic components} of  $T$. We denote by $\check T$ the acyclically indecomposable tournament indexing them.

 The next one is more involved:
\begin{theorem}  \label{morphology} There are twelve infinite acyclically indecomposable tournaments such that every infinite acyclically indecomposable tournament contains a  subtournament isomorphic to one of these tournaments.  
\end{theorem}

The
 twelve tournaments  mentionned in Theorem  \ref{morphology} are described in Section \ref{twelve}. At this point, we mention that they do not embed in each other, each one is countable and is the union of two acyclic tournaments. We also indicate that to an acyclic tournament  $C$ we associate  a set $\mathfrak B_{C}$ consisting of (at most) six tournaments denoted respectively   $C_{3 [C]}, V_{[C]}, T_{[C]}, U_{[C]}, H_{[C]}$ and  $K_{[C]}$. Let  $\mathfrak  B:= \mathfrak B_{\omega}\cup \mathfrak B_{\omega*}$ where $\omega$ is the tournament made of $\N$ and the natural (strict) order and $\omega^*$ is the dual of $\omega$. It turns out that each members of $\mathfrak  B$, except $K_{[\omega]}$, is acyclically prime and that $\check{K}_{[\omega]}$ is obtained from $K_{[\omega]}$ by identifying two vertices. The twelve tournaments mentionned in Theorem \ref{morphology} are obtained by  replacing  $K_{[\omega]}$ by $\check{K}_{[\omega]}$ in $\mathfrak B$. Indeed, as we will prove,  every infinite acyclically indecomposable  tournament
contains a member of $\mathfrak {B}$. The proof  is based on a separation lemma (Lemma \ref{separation})  and Ramsey Theorem.

Theorem  \ref{morphology} has a finitary version. Denote by $\underline n$ the acyclic tournament made of  $\{0, \dots, n-1\}$ with the natural (strict) order and set 
  $\check{\mathfrak B}_{\underline n}:=\{ \check{T}: T\in \mathfrak B_{\underline n }\}$ for each non negative integer $n$. 

\begin{theorem}\label{compactness}
For every non-negative integer $n$ there is an integer $a(n)$ such that every  finite tournament of size at least $a(n)$ which is acyclically indecomposable contains a subtournament isomorphic to a member of $\check{\mathfrak B}_{\underline n}$. 
\end{theorem}

An upper bound for $a(n)$ can be deduced from a careful analysis of the proof of Theorem \ref {morphology}. An existence proof is readily obtained  by means of the compactness theorem of first order logic. 

Indeed, suppose that the conclusion of Theorem \ref{compactness} is false.  Let $n$ be such that for every integer $m$ there is an acyclically indecomposable  tournament $T(m)$ of size at least $m$ which contains no subtournament isomorphic to a member of $\check {\mathfrak B}_{\underline n}$.  With the terminology  of \emph {embeddability}, we simply say that no member of $\check {\mathfrak B}_{\underline n}$ is embeddable in $T$. The compactness theorem of first order logic (or  an ultraproduct) yields a tournament $T$ such that for every first order-sentence $\varphi$ of the language of tournaments, $\varphi$  holds in $T$ whenever it holds in all of the $T(m)$, but finitely many. The fact that a given finite tournament is embeddable in a tournament can be expressed by the satisfaction of a first order sentence (in fact an existential one), thus  no member of $\check{\mathfrak B}_{\underline n}$   is embeddable into $T$. Our separation lemma ensures that the  fact that a tournament is acyclically indecomposable can be expressed by the satisfaction of a first order formula (Lemma \ref{logicindecomposable}). Hence $T$ is a acyclically indecomposable.  Since the size of the $T(m)$'s is unbounded, $T$  is infinite, thus,  from Theorem \ref{morphology}, some  $\check {X}_{[\alpha]}$ with $X_{[\alpha]}\in {\mathfrak  B}$ is embeddable in $T$. This tournament  is an increasing union of $\check {X}_{[\underline m']}$, for $m'\in \N$ (Corollary \ref{increasingunion}). Hence  $\check {X}_{[\underline n]} \in \check{\mathfrak {B}}_ {\underline n}$ is embeddable in $T$,  a contradiction. 

Let $\mathfrak A$ be the collection of tournaments  $T$ which can be written as a finite lexicographical sum of acyclic tournaments. Tournaments not in  $\mathfrak A$ are \emph{obstructions} to $\mathfrak A$. Clearly, no member of $\mathfrak A$ contains  a subtournament isomorphic to an obstruction. From Theorem \ref{morphology} (and the fact that $K_{[\omega]}$ can be embedded into  $\check{K}_{[\omega]}$)   $\mathfrak B$ is a  set of obstructions characterizing  $\mathfrak A$. And since the members of $\mathfrak B$ do not embed in each other, ${\mathfrak B}$ is a minimum sized set of obstructions. As a consequence: 
\begin{corollary} \label{morphology1}Let $T$ be an infinite  tournament,  then:
\begin{enumerate}[{-}]
\item Either $T$ is a lexicographical sum of acyclic tournaments indexed by a finite tournament.
\item Or $T$ contains as a subtournament a tournament isomorphic to a member of $\mathfrak B$.
\end{enumerate}
\end{corollary}

\subsection{Application to the profile of tournaments}
The {\it profile} of  a tournament $T$ is the function $\varphi_T$ which counts for each integer $n$ the number $\varphi_T(n)$ of  tournaments induced by $T$ on the  $n$-element subsets of $T$, isomorphic tournaments being identified.  The \emph{age} of  $T$ is the set $\mathcal A(T)$ of isomorphic types of subtournaments induced on the finite subsets of $V(T)$.   Clearly, the profile of $T$ depends only upon the age of $T$. 
We prove (see Section \ref {twelve} and Section \ref{section:computprofile}
):
\begin{lemma}\label {exponential} The ages of members of ${\mathfrak B}$ are six sets pairwise incomparable w.r.t. inclusion. For each $T\in {\mathfrak B}$  the growth of $\varphi_T$ is at least exponential, that is  $\varphi_T(n)\geq ac^n$ for some  reals $a>0$ and  $c>1$. 
\end{lemma}

It is easy to see that if $T$ is a lexicographical sum of acyclic tournaments indexed by a finite tournament, say $D$, then $\varphi_T$   is bounded from above by a polynomial (of degree at most $\vert D\vert -1$).  From Corollary \ref {morphology1} and Lemma \ref{exponential} we deduce:
\begin{theorem} \label{thetheorem}The profile of a tournament $T$ is either bounded from above by a polynomial, in which case  $T$ is a lexicographical sum of acyclic tournaments indexed by a finite tournament, or it growth is at least exponential.\end{theorem}

We give a precise description of the profile of a lexicographical sum of acyclic tournaments indexed by a finite tournament. 


 \begin{theorem}\label{main1} If  a tournament $T$ is a lexicographic sum of acyclic tournaments indexed by a finite tournament  then the generating series of the profile 
$$H_{\varphi_T}:= \sum_{n=0}^\infty \varphi_T(n)x^n$$ is a rational fraction of the form:
 \begin{displaymath}\frac{P(x)}{(1-x)(1-x^2)\cdots(1-x^k)}   \end{displaymath}
with $P(x)\in \Z[x]$ and $\varphi_T(n) \simeq an^{k-1}$ for some non-negative real $a$, the integer $k$ being the number of  infinite acyclic components of $T$.
\end{theorem}

The first part of Theorem \ref {main1} is a   consequence of a more general result about relational structures recently obtained by N.Thi\'ery and the second author \cite{pouzetthiery} that we record  in Section \ref{monodec}.

There are acyclically indecomposable tournaments of size $k$ for every integer $k\geq 3$, hence, according Theorem \ref{main1}, there are tournaments of arbitrarily large polynomial growth.

%


%






An other consequence of Theorem \ref{main1} is this:

\begin{corollary}
The growth of the profile of an infinite indecomposable  tournament is at least exponential. 
\end{corollary}

This research leaves open  the following:
\begin{problem} Find a  result,   similar to Theorem \ref{morphology},  for indecomposable tournaments and, possibly,  a finitary version. 
\end{problem}

The notion of acyclically indecomposable tournament was studied by J.F.Culus and B.Jouve in \cite{culus-jouve1}, \cite{culus-jouve2}, \cite{jouve}.  
The notion of profile was introduced in 1971 by the second author (see  \cite{fraisseclmt1}, \cite{fraissetr}) and developped in \cite{pouzet76},\cite{pouzettr}, \cite{pouzetrpe}; for a  survey  see \cite{pouzethammamet}.
The study of the orbital profile of  permutation groups is intensively studied by P.J.Cameron and his school  \cite {cameron1},  \cite{cameron2},  \cite{cameron3}.  The survey \cite{pouzethammamet} includes a presentation of Theorems \ref{morphology}, \ref{thetheorem} and  \ref{main1} with  an application to the structure of the age algebra of Cameron.  

This work was presented at the CGCS 2007(Luminy, France,  May 2-4 2007) in honor of Michel Deza. We are pleased to thanks Y.Manoussakis and the other  organisers of this conference for offering this opportunity to us. We thanks J.Nesetril for his comments and bibliographical references.

	This paper is  organized as follows. Section  \ref {perequisite} contains the material needed  about relational structures and tournaments.  Section \ref {section:acyclicdec} contains the main properties of acyclic decompositions of tournaments, particularly the proof of Theorem \ref{acyclicallyindecomposable1} and of our separation lemma.  Section \ref {proofTheoremmain1} contains the proof of Theorem  \ref{main1},  Section \ref{twelve} the description of $\mathfrak B_{\underline {n}}$ and $\mathfrak B$ with  their  main properties. Section {section:computprofile} contains the description of the profiles of members of $\mathfrak B$ and  Section \ref{sectiontheoremmorphology} the proof of Theorem \ref{morphology}.

\section{Perequisite}\label{perequisite}
We use standard set-theoretical notations. If $E$ is a set, $\vert E\vert$ denotes its cardinality. If $n$ is an integer, $[E]^n$ denotes the set of $n$-element subsets of $E$;  whereas $E^n$ denotes the set of $n$-tuples of elements of $E$. 

\subsection{Invariant structures and skew products}
A relation $\rho$ on a set $E$ is a map  from a finite  power $E^n$ of $E$ into the two element set $2:=\{0,1\}$; the integer $n$ is the {\it arity} of $\rho$, denoted $a(\rho)$ and $\rho$ is said $n$-ary.  If $n=2$ we say that $\rho$ is a {\it binary relation} and we denote $x \rho y$ the fact that $\rho(x,y)=1$. A {\it relational structure} is a pair  $ R:= (E, (\rho_i)_{i\in I})$ where each $\rho_i$ is a relation on  $E$. We denote by $R_{\restriction A}$ the relational structure induced by $R$ on $A$. A \emph{local automorphism} of $R$ is any isomorphism $h$ from an induced substructure of $R$ onto an other one. A pair $(E, \rho)$ where $\rho$ is a binary relation is a {\it directed graph}.  A {\it chain} is a pair $ C:= (A, \leq )$  where $\leq $ is a linear order on $A$. In this case  a local automorphism of  $C$ is 
every map $h$ from a subset $F$ of $C$ onto an other subset $F'$ of $C$ such that 
\begin{equation}\label{localaut}
x\leq y \Longleftrightarrow h(x)\leq h(y)
\end{equation}
for every $x,y \in F$. 
For each integer   $n$, let  $[C]^{n}$ be the set of  $n$-tuples $(c_{1},\dots,c_{n})$ of members of $A$ such that $c_{1}<\cdots <c_{n}$. These $n$-tuples will be identified with the $n$-element subsets of $A$.  If  $h$ is local automorphism of  $C$,  $F$ is its domain, $n$ is an integer  and $\vec{c}:= (c_1,\dots, c_n)\in [F]^n$, we will set  $\overline h(\vec {c}):= (h(c_1), \dots, h(c_n))$. 

Let  $\mathfrak L:= < C, R, \Phi>$ be a triple made of a chain  $C:= (A, \leq )$,  
a relational structure $ R:= (E, (\rho_i)_{i\in I})$ and a set  $\Phi$ of maps,   each one being a map $\varphi$  from  $[C]^{a(\varphi)}$ into  $E$, where $a(\varphi)$ is an integer, the  \textit{arity}  of  $\varphi$. 

We say that $\mathfrak L$ is  {\it invariant} if 
\begin{equation}
\rho_i(\varphi_1(\vec {a_1}), \dots, \varphi_m( \vec {a_m } ))= \rho_i(\varphi_1(\overline h(\vec{a_1})), \dots, \varphi_{m}(\overline h( \vec {a_{m} })))
\end{equation}
for every  $i\in I$, $m:=a(\rho_i)$, $\varphi_1,\dots, \varphi_m \in \Phi$, $\vec{a_1}\in [C]^{a(\varphi_1)}, \dots, \vec{a_m}\in [C]^{a(\varphi_m)}$ and every  local automorphism  $h$  of $C$ whose domain  contains the union of $\vec{a_1}, \dots, \vec {a_m}$.  

This technical condition expresses the fact that each $\rho_i$ is invariant under the transformations of the $a (\rho_i)$-tuples of $E$ which are  induced on $E$ by the local isomorphisms of $ C$.  
In the case of a  single binary relation $\rho$ and one $n$-ary function $\varphi$,  it says that $\varphi (\vec{a})\rho \varphi (\vec{b})$ depends only upon the relative positions of the components $\vec{a}$ and $\vec{b}$ on the chain $ C$.  

If $\mathfrak L:= <C,  R, \Phi>$ and $B$ is a subset of $A$, $\Phi_{\restriction B}:=\{\varphi_{\restriction [B]^{a(\varphi)}}: \varphi \in \Phi \}$ and  $\mathfrak {L}_{\restriction B}:= <C_{\restriction B},  R, \Phi_{\restriction B}>$  is  the {\it restriction} of $\mathfrak L$ to $B$. 

We will use the  following  straightforward consequence of Ramsey's theorem.\begin{theorem}\label{invar1}
Let $\mathfrak L:= < C,  R, \Phi>$ be a structure such that the domain $A$ of $ C$ is infinite, $ R$ consists of finitely many relations, and $\Phi$ is finite. Then there is an infinite subset  $A'$ of $A$  such that the structure 
$\mathfrak L_{\restriction A'}$ is   invariant.
\end{theorem}

Let $S: = (V, (\rho_{i})_{i \in I})$ be a relational structure and   $C: =
(A,
\leq)$ be a chain. A relational structure $R$
is a   {\it  skew product } of $S$ and $C$, denoted by $S
\bigotimes C$ if 
\begin{enumerate}
\item the domain is $A\times V$
\item for every  $x\in A $, the map  $v \rightarrow (x, v)$ is an isomorphism 
from  $S$ into $R$
\item for each local automorphism $h$ of $C$, the map $(h, 1_{V})$ defined by $(h,1_V)(x, v)=(h(x), v)$
is a local automorphism  of $R$.
\end{enumerate}
Let $\mathfrak L:= < C,  R, \Phi>$ where $\Phi:=\{\varphi_v: v\in V\}$ and $\varphi_v$ is the map from $A$ into $A\times V$ defined by $\varphi_v(x):= (x,v)$. Condition (2) expresses that $\mathfrak L$ is  
invariant. 

Theorem \ref{invar1} yields:
\begin{lemma} \label{lem:invar}Let $R$ be a  relational structure made of finitely many  relations and defined on a product $A\times V$. If $V$ is finite and $A$ is infinite,  then for every linear order $\leq $ on $A$ there is some infinite subset $A'$ of $A$ such that $R_{\restriction A'\times V}$ is a skew product of $R_{\restriction \{a\}\times V}$ and $C_{\restriction A'}$,   for some $a\in A$ and $C:= (A, \leq)$. 
\end{lemma} 

If $R$ is a skew product of a  finite relational structure $S$ by a  chain then $\varphi_R$ is bounded from above by some  exponential function. Indeed, if $S$ is such that  all its one-element restrictions are non-isomorphic, $\varphi(n)= \sum_{k=0}^{v}\varphi(n-k){n\choose k}$ where $v$ is the size of the domain of $S$, hence the generating series $\mathcal H_{\varphi_R}$ is a rational fraction and the result follows. If $S$ is arbitrary, its profile is dominated tems by terms by the previous one. It is not known if the  generating series of  a skew product $R$ of a  finite relational structure $S$ by a  chain is a rational fraction. It is not even known if the profile of $R$ is either polynomial or exponential. 

In this  paper,  we will consider skew product of a two-element tournament by a chain. For those which are not acyclic, their profile is asymptotically bounded from above by $\frac{1}{2}(1+\sqrt{2})^n$. As we  will see in Section \ref{section:computprofile}, their profile is bounded from below by some exponential.

 The notion of invariant structure appeared in \cite{charretton-pouzet}, Theorem \ref {invar1} was  an handy tool for using Ramsey's theorem. 
The notion of a skew product of a relation has appeared (under other names) in various papers of the second author (see \cite{pouzethammamet}). For some recent applications, see \cite{PSZ} and \cite{pouzet07}.

\subsection{Monomorphic decomposition of a relational structure}\label{monodec}
Let $R$ be a relational structure on $E$.  A
subset $B$ of $E$ is a \emph{monomorphic part} of $R$ if for every
integer $n$ and every pair $A, A' $ of $n$-element subsets of $E$ the
induced structures on $A$ and $A'$ are isomorphic whenever $A\setminus
B=A'\setminus B$. 
This notion has been introduced by N.Thi\'ery and the second author \cite{pouzetthiery}. We present the results we need. The following lemma gathers
the main properties of monomorphic parts.
\begin{lemma}\label{lemmahomogeneous}
  \begin{enumerate}[(i)]
  \item The empty set and the one element subsets of $E$ are
    monomorphic parts of $R$;
  \item If $B$ is a monomorphic part of $R$ then every subset of $B$
    too;
  \item Let $B$ and $B' $ be two monomorphic parts of $R$; if $B$ and
    $B'$ intersect, then $B\cup B'$ is a monomorphic part of $R$;
  \item Let $\mathcal B$ be a family of monomorphic parts of $R$; if
    $\mathcal B$ is up-directed (that is the union of two members of
    $\mathcal B$ is contained into a third one), then their union $B:=
    \bigcup \mathcal B$ is a monomorphic part of $R$.
  \end{enumerate}
\end{lemma}

\begin{corollary}For every $x\in E$, the set-union   $R(x)$  of all the monomorphic parts
of $R$ containing $x$ is a monomorphic part, the largest monomorphic part containing $x$.
\end{corollary}
\begin{proof}By $(i)$ of Lemma~\ref{lemmahomogeneous}, the set $R(x)$
contains $x$ and by $(iii)$ and $(iv)$ this is a monomorphic part,
thus the largest monomorphic part of $R$ containing $x$.
\end{proof}

We call the set $R(x)$  a \emph{monomorphic component of $R$}.

A {\it monomorphic decomposition } of  a relational structure $R$ is a partition $\mathcal P$ of $E$ into blocks such that for every
integer
$n$, the induced structures on two $n$-elements subsets $A$ and $A'$ of $E$ are isomorphic whenever the intersections $ A\cap B$  and
$ A'\cap B$ over each block $B$ of $\mathcal P$ have the same size.
\begin{proposition}\label{homogeneouscomponent}
 The monomorphic components of $R$  form a monomorphic decomposition
  of $R$ of which every monomorphic decomposition of
  $R$ is a refinement.
 \end{proposition}
We will call  {\it canonical} 
the decomposition of $R$ into monomorphic components.

Recently, N.Thi\'ery and the second author proved this: 

\begin{theorem}\label{theorem.quasipolynomial}
  Let $R$ be an infinite relational structure admitting  a monomorphic
  decomposition into finitely many blocks and let $k$ be the number of infinite blocks of the canonical decomposition of $R$, then:
  \begin{enumerate}
  \item The generating series
  $H_{\varphi_R}$ is a rational fraction of the form:
  \begin{displaymath}
    \frac{P(x)}{(1-x)(1-x^2)\cdots(1-x^k)} 
      \end{displaymath}
 where $P\in \Z[x]$.
 \item   $\varphi_R(n) \simeq an^{k-1}$ for some positive real $a$.
 \end{enumerate}\end{theorem}
  The proof of the first part and the fact that $\varphi_R(n) \simeq an^{k'}$ for some $k'\leq k-1$ is in \cite {pouzetthiery}. The proof that $k'=k-1$ was obtained afterward. 
  We will give below the proof of the second part  for the special case of tournaments. 

\subsection {Basic terminology and notations for tournaments} A {\it tournament} $T$ is a pair $(V, \mathcal E)$, where $\mathcal E$ is a binary relation on $V$ which is  irreflexive, antisymmetric and complete. Members of $V$ are the {\it vertices} of $T$, pairs $(x,y)$ of vertices such that $(x,y)\in \mathcal E$  are the {\it edges} of $T$. Given a pair   $u:= (x, y)$, resp. a set $\mathcal F$ of pairs,  we set  $u^{-1}:= (y,x)$, resp. $\mathcal F^{-1}:= \{u^{-1}: u\in \mathcal F\}$. The tournament $T^*:= (V, \mathcal E^{-1})$ is the \emph{ dual} of $T$.  If $A\subseteq V$, $T_{\restriction A}:= (A, \mathcal E \cap A\times A)$ is the tournament {\it induced on} $A$ or the {\it restriction of $T$ to $A$}. A {\it subtournament} of $T$ is any restriction of $T$ to a subset of $V$.   As usual, $V(T)$, resp. $E(T)$,  stands for the set of vertices, resp. edges,  of the tournament $T$. We also write $T(x,y)=1$ for $(x,y)\in E(T)$ and $T(x,y)=0$ for $(x,y)\not \in E(T)$.  An {\it isomorphism} from a tournament $T$ onto  a tournament $T'$  is a bijective map $f:V(T)\rightarrow V(T')$ such that $T(x,y)=T'(f(x),f(y))$ for all $(x,y)\in V(T)\times V(T)$. If $T'$ is a tournament, a   tournament $T$ is {\it isomorphic to} $T'$, resp. {\it is embeddable  into}  $T'$  if there is an isomorphism from $T$ onto $T'$, resp. onto a subtournament   of $T'$. A tournament is \emph{self-dual} if it is isomorphic to its dual.  A $3$-{\it element cycle},  or briefly a {\it $3$-cycle},  of a tournament $T$ is a $3$-element subset  $\{a,b,c\}$ of  $V(T)$ such that  $T(a,b)=T(b,c)=T(c,a)$.  The tournament induced on a $3$-cycle is also called a {\it $3$-cycle}. As a tournament, we will denote it by   $C_3$. A tournament $T$ is {\it acyclic} if no subtournament is a $3$-element cycle; this amounts to say that  the relation $E(T)\cup \{(x,x): x\in V(T)\}$  is a linear order. Up to reflexivity, acyclic tournaments and chains (alias totally ordered sets) being the same objects, we will use standard notions and notations used for chains.  For  example, we will say that the tournament  $(\N, <)$, where $<$ is  the strict order on $\N$,  has type  $\omega$; its dual is isomorphic to the tournament made of the set of negative integers equipped with the strict order, we will say that it has type   $\omega^*$.  Note that according to the theorem of Ramsey, every  infinite tournament contains a subtournament which is isomorphic to $\omega$ or to $\omega^*$. 

Let $D$ be a tournament and   $(T_i)_{i\in V(D)}$ be  a family of tournaments.  The {\it lexicographical sum  of the  tournaments $T_i$ indexed by the tournament $D$} is the tournament,  denoted $\sum_{i\in D}T_i$, and  defined as follows.   The vertex set is the disjoint union of the family $(V(T_i))_{i\in V(D)}$, that is $\bigcup \{V(T_i): i\in V(D)\}$ if the $V(T_i )$'s are pairwise disjoint and  $\bigcup\{ V(T_i)\times \{i\}: i\in V(D)\}$ otherwise.   Members  of this disjoint union being denoted by pairs $(x,i)$ with $x\in V(T_i)$, a pair $((x,i), (y,j))$ of vertices is an edge if either  $i=j$ and  $(x,y)\in E(T_i)$ or $(i,j)\in E(D)$. If $D$ has type $\omega$, resp. $\omega^*$, the lexicographical sum is an $\omega$-sum, resp. an $\omega^*$-sum.  If  $T_i =T$ for all $i\in V(D)$, this sum is a  {\it lexicographical product } of $T$ and $D$ denoted $T.D$. 

A subset $A\subseteq  V(T)$ of a tournament $T$ is {\it autonomous} if for every $x,x'\in A, y\not \in A$, $(x,y)\in E(T)$ if and only if $(x',y)\in E(T)$. The empty set, the one-element subsets and the whole vertex set  are autonomous and are said \emph{trivial}.  If $T$ has no other autonomous subset,   $T$ is {\it indecomposable} (an other denomination is  \emph{simple}, see \cite {erdos}, \cite {nesetril}). 
If $T$ is acyclic,  autonomous subsets coincide with intervals of the linear order, hence if $\vert V(T)\vert \geq 3$, $T$ is not indecomposable.  We also recall that if $T$ is a lexicographical sum  $\sum_{i\in D}T_i$, the subsets  of the form $V(T_i)$ are autonomous. Conversely, if the vertices of a tournament $T$ are partionned into autonomous subsets, then $T$ is the lexicographical sum of the  blocks of the partition. 
\newpage
\section {Acyclic decompositions  of tournaments}\label{section:acyclicdec}
\subsection{Proof of Theorem \ref{acyclicallyindecomposable1}.} We  recall the following result (which holds for arbitrary binary relations).

\begin{lemma} \label{autonomous}Let $T$ be a tournament. 
\begin{enumerate}
\item The union of two autonomous subsets of $T$ with a non empty intersection is autonomous.
\item The union of a family $\mathcal F$ of automous subsets of $T$ which is  closed under finite union  is autonomous. 
\end{enumerate}
\end{lemma}

\begin{lemma}\label{union}
The union of two acyclic autonomous subsets of a tournament  is acyclic. 
\end{lemma}
 The proof is immediate.

Applying Lemma \ref{autonomous} and  Lemma \ref{union},  we get:
 \begin{lemma}\label{keydec} Let $T$ be a tournament and $x\in V(T)$. Then  the set-union $Ac(T)(x)$ of all the acyclic autonomous subsets of $T$ containing $x$ is the largest  acyclic autonomous subset containing $x$. 
 \end{lemma}

An {\it acyclic component} of  $T$  is any subset of $V(T)$ of the form $Ac(T)(x)$. 

From Lemma \ref{keydec}, we also have immediately: 
\begin{lemma} \label{acyclic1}
Let $T$ be a tournament. Then:
\begin{enumerate}
\item Every acyclic autonomous subset is contained into an acyclic  component.
\item  The  acyclic components of $T$ form a partition of $V(T)$ into autonomous subsets . 
\item  Every partition of $V(T)$ into acyclic autonomous subsets is a refinement of the partition into acyclic components.  \end{enumerate}
 \end{lemma}

As a corollary of Item (3) we get:

\begin{proposition}\label{acyclic} A tournament is a lexicographical sum of acyclic tournaments indexed by a finite tournament if and only if it has only finitely many acyclic components. 
\end{proposition}
Let  $Ac(T)$ be the set of acyclic components of a tournament $T$; set $ac(T):= \{\vert A\vert : A\in Ac(T)\}$ and $\overline {ac}(T)$ be the sequence of the  elements of $ac(T)$ sorted in a decreasing order.  As a corollary of the existence of acyclic components, we get: 

\begin{corollary}\label{spectrum}
If $T$  and $T'$ are two isomorphic  finite tournaments,  $\overline {ac}(T)=\overline{ac}(T')$. 
\end{corollary}

  Let $T$ be a tournament and let  $p: V(T)\rightarrow Ac(T)$ defined by setting $p(x):= Ac(T)(x)$.  Let $\mathcal E:= \{(Ac(T)(x), Ac(T(y)):  Ac(T)(x)\not =Ac(T(y) \; \text {and} \; T(x,y)=1\}$ and let $\check{T}:= (Ac(T), \mathcal E)$.

\begin{lemma} 
Let $T$ be a tournament. Then $\check T$ is acyclically indecomposable and   $T$ is the lexicographical sum of its acyclic components indexed by $\check T$.
\end{lemma} 
\begin{proof} According to Item (2) of Lemma  \ref{acyclic1} above,  $T$ is  the lexicographical sum of its acyclic components indexed by $\check{T}$.  Let us prove that   $\check {T}$ is acyclically indecomposable. Let $A$ be an acyclic autonomous subset of $V(\check{T}):=Ac(T)$. Then $\bigcup A$ is an autonomous subset  of $T$. Since  $T_{\restriction \bigcup A} $ is a lexicographic sum of acyclic tournaments indexed by the acyclic tournament $\check{T}_{\restriction A}$, $\bigcup A$ is also acyclic. Consequently, $\bigcup A$ reduces to a single acyclic component and $A$ is a singleton. This proves our assertion. 
\end{proof}

With this lemma the proof  of Theorem  \ref{acyclicallyindecomposable1} follows. 
\subsection{A separation lemma}
A {\it diamond}, resp. a   {\it double diamond},  is a tournament  obtained by replacing a vertex of a $2$-element tournament, resp. a $3$-element acyclic tournament,   by a $3$-element cycle.  A   double diamond is  self-dual  if and only if  the middle element of the $3$-element acyclic tournament is replaced by a $3$-element cycle. 

We have the following separation lemma which generalizes Lemma \ref{union}:
\begin{lemma}\label{separation}Two vertices $x,y$  of  a tournament $T$ do not belong to an  acyclic autonomous subset   of $T$ if and only if $x$ and $y$ are distinct and either: 
\begin{enumerate}
\item $x$ and $y$  belong to some $3$-element cycle, or
\item  $x$ and $y$ belong to some diamond, or 
\item $x$ and $y$ belong to some self-dual double diamond.  
 \end{enumerate}
\end{lemma}
\begin{proof} Let $x, y$ be two distinct vertices of $T$  and $A$ be an autonomous acyclic subset of $V(T)$ containing $x,y$. If $x,y$ belong to some $3$-element cycle $C$, then, since $A$  is acyclic, the 
third element, say $z$, does not belong to $A$. Since C is a cycle, $T (z, x)= T (z, y)$
whereas, since $A$ is autonomous, $T (z, x) = T (z, y)$. A contradiction. 
 If  $x,y$ belong to some diamond $\delta$, then from the previous case, they cannot belong to the $3$-cycle of the diamond. With no loss of generality we may suppose that $\delta$ is a \emph{positive} diamond, that is  $\delta= \delta ^+:= (\{a,b,c, d\}, \{(a,b), (b,c), (c,a), (a,d), (b,d),(c,d)\})$, with $x:=a, y:= d$. 
Since $A$ is autonomous, $a,d\in A$ and $T(b,a)\not =T(b,d)$,  we get  $b\in A$.  Hence, $A$ contains two vertices of the $3$-cycle $\{a,b,c\}$. From the previous case, this contradicts the fact that $A$ is acyclic.  If $x, y$ belong to a self-dual double diamond $D:= \{a,b,c, d, d'\}$, then  either they belong to one of  the two diamonds included in $D$ or they coincide with the end-points $d, d'$ of $D$.  As seen above, the first case contradicts the fact that $A$ is acyclic. In the second case, since $A$ is autonomous and $T(z, d)\not =T(z,d')$ for every element $z$ of the $3$-cycle, $A$  must contain each element of the $3$-element cycle, hence it cannot be acyclic.

Conversely,  suppose that $x$ and $y$ does not belong to an acyclic autonomous subset of $T$, in particular $x \not =y$. 
Let $Z:=\{x,y\}$, $Z_i:= \{z\in V(T)\setminus Z: T(z,x)=T(z,y)=i\}$ for $i\in \{0, 1\}$ and $Z_{\frac{1}{2}}:=  \{z\in V(T)\setminus Z: T(z,x)\not =T(z,y)\}$. With no loss of generality, we may suppose $T(x,y)=1$. 

\noindent \begin{claim} \label{claim:claim1}  
\begin{enumerate}[{(a)}]
\item $Z_{\frac{1}{2} }\not =\emptyset$.
\item If $x,y$ does not belong to a diamond or a $3$-cycle then   $Z_{\frac{1}{2}}$ is an autonomous subset of  $Z\cup Z_{\frac{1}{2}}$ and  $Z\cup Z_{\frac{1}{2}}$ is  an autonomous subset of $T$. 
\item Furthermore, $Z_{\frac{1}{2}}$ contains a $3$-cycle and $x$, $y$ belong to a double diamond. 
\end{enumerate}
\end{claim}
{\bf Proof of Claim \ref{claim:claim1}.}

(a) Since $Z$ is acyclic, it cannot be autonomous. Hence,  there is $z\not \in Z$ such that $T(z, x)\not = T(z, y)$. We have  $z\in Z_{\frac{1}{2} }$, hence  $Z_{\frac{1}{2} }\not =\emptyset$.

(b) Suppose that $x$ and $y$ does not belong  to a $3$-cycle. Let $z\in Z_{\frac{1}{2} }$. Since $T(z, x)\not =T(z,y)$ we have $T(x, z)=T(z,y)$.  Since $\{x,y,z\}$ is not a $3$-cycle,  we have $T(x,z)=T(x,y) =T(z, y)$ . Since the values $T(x,z)$ and $T(z,y)$ do not depend upon our choice of $z$, $Z_{\frac{1}{2}}$ is an autonomous subset of  $Z\cup Z_{\frac{1}{2}}$. Let $i\in \{0, 1\}$ and  $z_i\in Z_i$ . If $T(z_i, z)\not = i$ for some $z\in Z_{\frac{1}{2} }$, then $\{x,y,z, z_i\}$ is a diamond. Hence, supposing that  $x,y$ does not belong to a diamond, we have $T(z_i, z)=i$ for all $z \in Z\cup Z_{\frac{1}{2}}$, proving that $Z\cup Z_{\frac{1}{2}}$ is autonomous.  

(c) From our hypothesis, $Z\cup Z_{\frac{1}{2}}$ cannot be acyclic. Since $Z_{\frac{1}{2}}$ is an autonomous subset of  $Z\cup Z_{\frac{1}{2}}$, $Z_{\frac{1}{2}}$ and no  $3$-cycle contains $\{x,y\}$, $Z_{\frac{1}{2}}$ cannot be acyclic.  Let $C$ be  a $3$-cycle included into $Z_{\frac{1}{2}}$. Then  $C\cup Z$ is double diamond containing $x$ and $y$, as claimed.
\endproof

With this claim, the proof is complete. 
\end{proof}
 \begin{lemma}
Let $T$ be a tournament , $A$ be  a subset of $V(T)$ and $\kappa$ the number of acyclic components of $T$ which meet  $A$. Then there is a subset $A'$ of $V(T)$ containing $A$ such that $\vert A'\setminus A\vert \leq 3.{{\kappa}\choose {2}}$ and the acyclic  components of $T_{\restriction A'}$  are the traces on $A'$ of the acyclic components of $T$. \end{lemma}

\begin{proof} Let $U:= \{X\in Ac(T): X\cap A\not =\emptyset\}$. For each $X\in U$, select $a_X\in X\cap X$. According to Lemma \ref{separation}, for each pair of distinct elements $a_X$, $a_Y$,  we may select a subset $F_{X,Y}$ containing $a_X$ and $a_Y$ such that $T_{\restriction F_{X,Y}} $ is either a $3$-element cycle, a diamond, or a self-dual double diamond. Set  $A':= A\cup \bigcup\{F_{X,Y}: \{X,Y\}\in [U]^2\}$ and $T':= T_{\restriction A'}$.  The traces on $A'$ of the acyclic decomposition of $T$ form a partition of $A'$ into acyclic autonomous subsets. \end{proof}

\begin{corollary}If $T$ is an acyclically indecomposable tournament, every subset $A$ of $V(T)$ extends to a subset $A'$ such that $T_{\restriction A'}$ is acyclically indecomposable and $\vert A'\setminus A \vert \leq 3. {\vert A\vert  \choose 2}$.
\end{corollary}

 \begin{corollary} \label{traces}Let $T$ be a tournament and $A$ be a subset of V(T). If $A$ meet each acyclic component of $T$, then the acyclic components of $T_{\restriction A}$ are the traces on $A$ of the acyclic components of $T$. 
\end{corollary}

\subsection{Relation with logic formulas}

Let $L$ be the first order language with equality  for which the only non logical symbol is a binary predicate denoted $<$.  In this language, tournaments are models of a universal sentence, namely the sentence $\theta:= \forall x \forall y (((x=y)\vee (y<x)\vee(x<y))\wedge ((\neg y<x)\vee(\neg x<y)))$

\begin{lemma} There is a two-variables first-order formula $\phi(x,y)$ of the language of tournaments such that for every tournament $T$ and  every $(a, b)\in V(T)\times V(T)$, the pair $\{a, b\}$ is no included into an acyclic autonomous subset of $T$ if and only if $T$ satisfies $\phi(a,b)$. 
\end{lemma}
\begin{proof}
Set $\phi(x,y):=\phi_1(x,y)\vee \phi_2(x,y)\vee \phi_3(x,y)$ such that the satisfaction of $\phi_1(a,b)$, resp. $\phi_2(a,b)$, resp. $\phi_3(a,b)$,  in a tournament $T$,  expresses that $a$ and $b$ are two vertices of a $3$-cycle, resp. a diamond, resp. are the end-vertices of a self-dual double diamond.  For an example, $\phi_1(x,y):=\theta(x,y) \vee \theta(x,y)$, where $\theta (x,y):= x<y \wedge (\exists z (y<z \wedge z<x))$.  This extends easily to $\phi_2(x,y)$ and $\phi_3(x,y)$. In particular,  $\phi(x,y)$ is a universal sentence.
\end{proof}

 \begin{lemma} \label{logicindecomposable}There is a first-order sentence $\phi$ of the language of tournaments such that a tournament $T$ is acyclically indecomposable if and only if it satisfies $\phi$.
\end{lemma}
\begin{proof}
Set $\phi:=\forall x\forall y \phi(x,y)$. 
\end{proof}

\subsection{Acyclic components and monomorphic parts}

Let $a, b$ be two distinct vertices of a tournament $T$. Let $C(a, b)$ be the set of vertices $x$ such that $\{a,b,x\}$ is a $3$-cycle of $T$ (that is $C(a,b) := \{x: T(x, a)=T(a,b)=T(b,x)\}$. 

\begin{lemma}\label{monoacypart}Let $T$ be a tournament and $A$ be a  subset of $V(T)$. 
\begin{enumerate}
\item If $A$ is acyclic and autonomous then $A$ is a monomorphic part  of $T$.
\item  If $A$ is a monomorphic part and no pair of distinct  vertices  of $A$  belongs to a $3$-cycle of $T$  then $A$ is included into an acyclic component of $T$. 
\item  If $A$ is a monomorphic part, there is a $3$-cycle which contains two vertices of $A$ if and only if $A$ is included into some  $3$-cycle  of $T$ and 
\begin{enumerate}[{-}]
\item either $A$ is an autonomous $3$-cycle,
\item or $A=\{a,b\}$ ,  $C(a,b)$ is acyclic and $\{a,b\}\cup C(a,b)$ is autonomous in $T$ . 
\end{enumerate}
\item If $A$ is a monomorphic component then either $A$ is an autonomous $3$-cycle or $A=\{x,y\}$ and $A\cup C(\{a,b\})$ is autonomous, or $A$ is an acyclic component of $T$.
\end{enumerate}
\end{lemma} 
\begin{proof} 
Assertion (1)  is obvious. 

Assertion (2).  If $A$ was no included into an acyclic component of $T$  then,  according to Lemma \ref{separation},  two  distinct  vertices $x$ and $y$  of $A$ would  belong  either to a  $3$-cycle of $T$, or to some diamond, or to 
some self-dual double diamond.  The first case is excluded by our hypothesis. The two other cases cannot happen. Indeed, if $x$ and $y$ belong to some diamond  $\delta$, then they cannot belong to the $3$-cycle of the diamond. With no loss of generality we may suppose that $\delta$ is a positive diamond, eg $\delta= \delta ^+:= (\{a,b,c, d\}, \{(a,b), (b,c), (c,a), (a,d), (b,d),(c,d)\})$, with $x:=a, y:= d$. 
Since $A$ is a monomorphic part, the tournaments  $T_{\restriction \{a, b, c\}}$ and $T_{\restriction \{d,b, c\}} $ must be isomorphic, which is impossible since the first one is a $3$-cycle and the other is acyclic. Thus this case cannot happen. If $x$ and $y$ belong to some double diamond $D:=\{ab,c,d,d'\}$, then since the previous cases cannot happen, $x$ and $y$ are the end-poinds $d$, $d'$  of the double diamond.  The two tournament $T_{\restriction \{a,b,c,d\}}$ and $T_{\restriction \{a,b,c, d'\}}$ are opposite, hence cannot be isomorphic. Whereas, since  $A$  is a monomorphic part, they must be isomorphic. A contradiction.

Assertion (3). Let  $A$ be  a monomorphic part of $T$. Let $C$ be a $3$-cycle containing at least two vertices $\{a, b\}$  of $A$ and let $c$ be the third vertex of $C$. 

\begin{claim}\label {claim:monoacypart} $A\subseteq C$. 
\end{claim}
Indeed, otherwise let $y\in A\setminus C$.  Since $C$ contains, at least,  two elements 
of $C$, there is one, say $a$,  such that $T(c, a)=T(c,y)$. Hence $D:=\{a, c,y\} $ is an acyclic tournament. Since  $D\setminus A=C\setminus A$ and $A$ is a monomorphic part, $D$ and $C$ must be isomorphic, which is not the case. This proves our claim. 

{\bf Case 1.} $A$ has three elements.  From Claim \ref{claim:monoacypart} above,  $C=A$.  Let  $y\not \in A$. Since $\vert A\vert =3$, there are   two element $a, b \in A$ such that $T(a,y)=T(b, y)$. We claim that for the remaining   element $c$ of $A$, $T(c, y)=T(a, y)$. Otherwise,  $T(c, y)=T(y, a)$ and since from the claim above,   $\{ a, c,  y\}$ cannot be a $3$-cycle, $T(y,a)=T(c, a)$. By the same token, $T(c, b)=T(y, b)$, hence $T(c,a)=T(c, b)$, contradicting the fact that $A$ forms a $3$-cycle. This proves our claim.  It follows that $A$ is autonomous. 

{\bf Case 2}. $A$ has two elements. First $C(a, b)$ is acyclic. Otherwise, if $D$ be a $3$-cycle included into $C(a,b)$ then $\{a\}\cup D$ and $\{b\}\cup D$ are two opposite diamonds, hence cdannot be isomorphic, contradicting the fact that $A$ is a monomorphic part. Next, $X:= A\cup C(a,b)$ is autonomous. Let $x\in V(T)\setminus X$ and $z\in C(a,b)$. We claim that the tournament $T_{\restriction \{a,b,x,z\}}$ is a diamond. Indeed, as a $4$-vertices tournaments, it contains at most two $3$-cycles. Since it contains the $3$-cycle  $\{a, b,z\}$, the restrictions  $T_{\restriction \{a, z,x\}}$ and $T_{\restriction  \{b, x,z\}}$- which are isomorphic since $A$ is a monomorphic part- must be acyclic. Since $x\not \in C(a,b)$, $\{a,b,x\}$ cannot be a cycle. Since $T_{\restriction \{a,b,x,z\}}$ contains just one $3$-cycle, this is a diamond. hence $T(x,a)=T(x,b)=T(x,z)$, proving that $X$ is autonomous.

Conversely, one can easily check that if $A$ satisfies the stated conditions then it is a monomorphic part. 

Assertion (4) follows immediately from Assertion (2). 

\end{proof}

From Assertion (1) and Assertion (4) we get:

\begin{corollary}\label{acyclicmonomorphic}
Let $T$ be a tournament and $A$ be a  subset of $V(T)$. If $\vert A\vert \geq 4$ then $A$ is an acyclic component of $T$ if and only if $A$ is a monomorphic component of $T$. 
\end{corollary}

\section{Proof of Theorem \ref{main1}} \label{proofTheoremmain1} 

Let $T$ be a tournament which is a lexicographical sum of finitely many  acyclic tournaments. Let $p$ be the number of acyclic components, $k$ be the number of the  infinite components. According to Lemma \ref{monoacypart}, each acyclic component is a monomorphic part  of $T$, and according to Corollary \ref{acyclicmonomorphic}, $T$ has exactly $k$ infinite monomorphic components, 
hence from part 1 of Theorem \ref{main1}, the generating series 
  $H_{\varphi_T}$ is a rational fraction of the form:
  \begin{displaymath}
    \frac{P(x)}{(1-x)(1-x^2)\cdots(1-x^k)} 
  \end{displaymath}
 where $P\in \Z[x]$.
 Furthermore, $\varphi_R(n) \simeq an^{k'}$ for some $k'\leq k-1$. 
 To complete the proof of Theorem  \ref{main1}, it remains to prove that $k'=k-1$. This is a consequence of Proposition \ref{growthk} below.

Let $n$ be a positive integer. A {\it partition of $n$} is a finite decreasing sequence $x_1\geq\cdots \geq x_k$ of positive  integers such that $x_1+\cdots +x_k=n$. The integers in this sequence are the {\it parts of the partition}. Set $\mathfrak p_{k}(n)$ for the number of partitions of the integer 
$n$ into at most $k$ parts, and set $\mathfrak p_{k}(0):=1$. As it is well-known, the generating series ${\mathcal H}_{\mathfrak {p}_{k}}:= \sum_{n=0}^{\infty}  \mathfrak p_{k}(n)x^{n}$ is the rational fraction $\frac {1}{(1-x)\cdots (1-x^{k})}$ and 
 $\mathfrak p_{k}(n) \simeq
\frac{n^{k-1}} {(k-1)!k!}$. We also recall that  {the partition function $\mathfrak p$}  counts  the number $\mathfrak p(n)$ of partitions of the integer 
$n$. A famous result of Hardy and Ramanujan,  1918, asserts that
$\mathfrak p(n) \simeq
\frac{1}{4
\sqrt{3n}} e^{\pi \sqrt{\frac{2n}{3}}}$. 

\begin{proposition}\label{growthk} If a   tournament $T$ is a finite lexicographical sum of acyclic tournaments, then 
\begin{equation}\label{partition}
\varphi_T(n)\geq \mathfrak {p}_k(n-p)
\end{equation} where $p$ is  the number of acyclic components of $T$, $k$ is the number of the infinite one and $n\geq p$. 
In particular the growth of $\varphi_T$ is at least a polynomial with degree $k-1$. 
\end{proposition}

\begin{proof}
 Let $Ac(T):= \{A_1, \cdots, A_p\}$ be the set of acyclic  components of $T$, enumerated in such a way that $A_1,\dots, A_k$ are infinite. Let $n\geq p$. To a  decreasing sequence $\vec{x}:= x_1\geq\cdots \geq x_{k'}$ of positive  integers such that $x_1+\cdots +x_{k'}=n-p$ and $k'\leq k$ associate the $p$-element sequence $1+\vec{x} = (1+x_1, \dots, 1+x_k, 1, \dots, 1)$ and a subset $A_{\vec{x}}$ of $V(T)$ such that $\vert A_{\vec{x}}\cap A_i\vert =x_i+1$ if $i\leq k'$ and $\vert A_{\vec{x}}\cap A_i\vert =1$ otherwise. Set $T_{\vec{x}}:=T_{\restriction A_{\vec{x}}}$. 

{\bf Claim} If $\vec{x}\not =\vec{x'}$ then $T_{\vec{x}}$ and $T_{\vec {x'}}$ are not isomorphic. 

Since $T$  contains at least an element of each acyclic component of $T$, the acyclic decomposition of $T_{\vec{x}}$ is induced by the acyclic decomposition of $T$ (Corollary \ref{traces}). Hence, $\overline {ac}(T_{\vec{x}})= 1+\vec{x}$. If  $T_{\vec{x}}$ and $T_{\vec{x'}}$ are isomorphic, $ac(T_{\vec{x}})=ac(T_{\vec{x'}})$ (Corollary \ref{spectrum}), thus  
 $\overline {ac}(T_{\vec{x}})=\overline {ac}(T_{\vec{x'}})$, that is $1+\vec{x}=1+\vec{x'}$ which yields $\vec{x}=\vec{x'}$.  
 
 Inequality (\ref{partition}) follows immediately.\end{proof}

\section{Twelve tournaments}\label{twelve}

Let  $C:= (A, <)$ be an acyclic tournament. We define six tournaments, denoted respectively by $C_{3 {[C]}}$, $V_{[C]}$, $T_{[C]}$, $U_{[C]}$,  $H_{[C]}$ and  $K_{[C]}$.

$\bullet$ The tournament $C_{3 {[C]}}$ is the lexicographical product  $C_3. C$ of $C_3$ by  $C$. 

$\bullet$ The vertex set of  $V_{[C]}$ is $A\times  2 \cup \{a\}$, where ${2}:=\{0, 1\}$ and  $a\not \in A\times  2$. A pair $(e,e')$ of vertices is an edge of $V_{[C]}$ in the following cases:
\begin{enumerate}  [{(i)}]
\item $e=(x,i)$, $e'= (x',i')$ and  either   $x<x'$ or $x=x'$ and  $i<i'$; 
 \item  $e= a$, $e'=(x',0)$; \item  $e=(x,1)$ and $e'=a$.  
 \end{enumerate}


$\bullet$ The four remaining tournaments  have the same vertex set, namely $A\times 2$. In order to define their  edge sets, let  $i\in  {2}$. Set $h_{i}:= ((0,i), (1,i))$, $v_{i}:= ((i,0), (i,1))$,  $d_{i}:= ((0,i),(1,i+1))$ (where $i+1=1$ if $i=0$ and $0$ otherwise). Let $X \subseteq  2\times  2\setminus \{ v_1, v_1^{-1}\}$. Let  $\Delta(C,X)$ be the directed graph whose vertex set is $A\times  2$ and edge set the union of the following three  sets:
\begin{enumerate} [{(a)}]
\item $\{((x,i), (x,j)):  ((0,i),(0,j))\in X\}$;
  
\item $\{((x,i), (y,j)): (x,y)\in E(C) \text{ and}   ((0,i),(1,j))\in X\}$;

\item $\{((x,i), (y,j)): (y,x)\in E(C) \text{ and}   ((1,i),(0,j))\in X\}$.
\end{enumerate}
Set $Y:= \{h_{0}, v_{0}\}$.   If $X:=\{d_{0}^{-1},  d_{1}^{-1}, h_{1}\} \cup Y$, resp. $X:= \{d_{0}^{1}, d_{1},  h_{1}^{-1}\}\cup Y$, resp.  $X:= \{ d_{0}^{-1}, d_{1}, h_{1}\}\cup Y$, resp. $X:= \{ d_{0}^{-1}, d_{1},
h_{1}^{-1}\}\cup Y$ then $\Delta(C,X)$ is a tournament  denoted by  $T_{[C]}$,  resp. $U_{[C]}$,  resp. $H_{[C]}$, resp. $K_{[C]}$. 

Conditions $(a),(b), (c)$ above simply mean that $\Delta (C,X)$ is a skew product of a binary relation on $\{0,1\}$ by the chain $(A, \leq)$. We choose  $X$ in such a way that $\Delta (C,X)$ is a tournament and in fact a skew product of the tournament   $\underline 2$ (for which $0<1$) by $(A, \leq)$. Deciding furthermore that this tournament will contains all pairs $((x,0), (y,0))$ such that $x<y$, we have only eight possible choices for the three remaining   pairs belonging to $X$. It  turns out that three choices yield acyclic tournaments. On the remaining five choices, two tournaments are dual of each other, namely $U_{[C] }$ and $U'_{[C]}:=\Delta (C,X)$ where $X:= \{ d_{0}^{-1}, d_{1}^{-1},
h_{1}^{-1}\}\cup Y$. We do  not  need to add to our list  tournaments of the form $U'_{[C]}$. Indeed, our aim is to obtain a minimal list of unavoidable infinite acyclically indecomposable tournaments. And  it follows from Item (\ref{item:2CleqC}) of Lemma \ref{lem:basictournament} below,  that $U_{[\omega^*]}$ is embeddable $U'_{[\omega]}$ and $U_{[\omega]}$ is embeddable $U'[\omega^*]$. 
 
\begin{lemma}\label{lem:basictournament}Let $C$ be an acyclic tournament, then:
(i) $(C_3.C)^*$ is isomorphic to $C_3. (C^*)$;  
(ii)  $(V_{[C]})^*$ is isomorphic to $V_{[C^*]}$; 
(iii) $(T_{[C]})^*$ is isomorphic to $T_{[C^*]}$; 
(iv) \label {item:2CleqC} If $\underline 2.C$ is embeddable in $C$ then $(U_{[C] })^*$ and $U_{[C^{*}]}$ are embeddable  in each other; (v) $(H_{[C]})^*$ is isomorphic to $H_{[C^*]}$;
 (vi) $K_{[C]}$ is   self dual.
\end{lemma}
\begin{proof}
We only check Assertion (iv).
We claim that  $(U_{[C]})^{*}$ is embeddable  in $U_{[C^{*}]}$. Indeed, let $\varphi:A\times 2\rightarrow (A\times 2)\times 2$ defined by $\varphi((x,i)):= ((x,i), i)$.  Then, as it is easy to check, $\varphi$ is an embedding from $(U_{[C]})^{*}$ into $U_{[(2.C)^*]}$. From our hypothesis $(\underline 2.C)^*$ is embeddable in $C^*$,  thus $U_{[(2.C)^*]}$ is embeddable in $U_{[C^*]}$. This proves our claim. Applying this claim to $C^*$ we get that $U_{[C^{*}]}$ is embeddable in $(U_{[C]})^{*}$ as required. 
\end{proof}

 We denote respectively by $\mathfrak C_{3}$, $\mathfrak V$, $\mathfrak T$,  $\mathfrak U$, $\mathfrak H$,  and $\mathfrak K$ the  collections of tournaments  $C_{3 {[C]}}$, $V_{[C]}$, $T_{[C]}$,  $U_{[C]}$, $H_{[C]}$  and  $K_{[C]}$ when $C$ describe all possible acyclic tournaments. We denote by   $\mathfrak C_{3,  <\omega}$, resp. $\mathfrak V_{<\omega}$,  $\mathfrak T_{<\omega}$, $\mathfrak{U}_{<\omega}$,  $\mathfrak H_{<\omega}$, $\mathfrak K_{<\omega}$,  the collection of finite tournaments which are embeddable  into some member of the corresponding collection.

 Some members of $\mathfrak V_{<\omega}$, $\mathfrak T_{<\omega}$ and  
 $\mathfrak U_{<\omega}$ and  have been considered previously. We will refer to  some known properties of these tournaments. We use the presentation given in \cite {boudabbous}.
Let  $h \geq 2$ be an integer, denote by 
$T_{2h+1}$, $U_{2h+1}$ and $V_{2h+1}$  the tournaments defined on
$\{0,\ldots, 2h\}$  as follows.\\ (i)$T_{2h+1 \restriction \{0,\ldots,h\}} =
U_{2h+1 \restriction \{0,\ldots,h\}} = 0<\ldots<h$,
$T_{2h+1\restriction \{h+1,\ldots,2h\}}=
(U_{2h+1})^{\ast}_{\restriction \{h+1,\ldots,2h\}} = h+1<\ldots<2h$.\\
(ii) For every
$i \in \{0,\ldots,h-1\}$, if $j \in \{i+1,\ldots,h\}$ and if $k
\in \{0,\ldots,i\}$, then $(j,i+h+1)$ and $(i+h+1,k)$
belong to $E(T_{2h+1})$ and $E(U_{2h+1})$.\\
(iii)$V_{2h+1 \restriction \{0,\ldots,2h-1\}} = 0<\ldots<2h-1$ and for $i \in
\{0,\ldots,h-1\}$, $(2i+1,2h)$ and $(2h,2i)$ belong to
$E(V_{2h+1})$.

According to Schmerl and Trotter  \cite{schmerl}, these tournaments  are
indecomposable and moreover  a finite tournament $T$  on at least five vertex is \emph{critically indecomposable} (in the sense that $T$ is indecomposable and for every $x \in V(T)$, the subtournament $T_{\restriction V(T)\setminus \{x\}}$ is not indecomposable) if and only if it is isomorphic to one of these tournaments.

We will need the following result \cite{boudabbous}.
\begin{lemma} \label{lem:TUV}Given three integers
$h_{1}$, $h_{2}$, $h_{3} \geq 2$, the tournaments $V_{2h_{3}+1}$, $T_{2h_{1}+1}$ and
$U_{2h_{2}+1}$  are
incomparable with respect to embeddability.
\end{lemma}

These tournaments belong to  $\mathfrak V_{<\omega}$, $\mathfrak T_{<\omega}$ and $\mathfrak U_{<\omega}$. Indeed: 

\begin{fact} \label{fact:VTU}$V_{2h+1}$ is isomorphic to $V_{\underline h}$, 
$T_{2h+1}$ is isomorphic to $T_{[\underline{h+1}]}$ minus the vertex $(0, 1)$ and $U_{2h+1}$ is isomorphic to $U_{[\underline{h+1}]}$ minus the vertex $(0, 0)$. \end{fact}
\begin{lemma}\label{lem:twelveincec} (i)\label{item:c3}   If $C:=(A, <)$ is an non-empty acyclic tournament,  $C_3[C]$ is acyclically indecomposable and not indecomposable except if $\vert A \vert= 1$. In fact, no indecomposable subset of $C_3[C]$ has more than three elements. (ii)$V[C]$ is indecomposable, hence acyclically indecomposable. 
(iii)$T_{[C]}$ is indecomposable, except if $\vert A\vert \geq 2$ and $C$ has a least and largest element. In this latter case $\{(m,0), (M,1)\}$  (where $m$ and $M$ are  the least and largest element of $C$) is an acyclic component, $T_{[C]}$ minus the vertex $(m,0)$ is isomorphic to $\check{T}_{[C]}$ and is indecomposable. (iv) $U_{[C] }$ is acyclically indecomposable for $\vert A\vert \geq 2$. If moreover  $C$ has a least element  $U_{[C] }$  is not indecomposable, $U_{[C] }$ minus  the vertex $(m,0)$ (where $m$ is the least element of $C$) is isomorphic to $\check{U}_{[C]}$ and  is indecomposable. (v)\label{item:h}  $H_{[C]}$ is indecomposable except for $\vert A\vert =2$.
(vi) $K_{[C]}$ is  never indecomposable; in fact its indecomposable subsets have at most three elements. It is acyclically indecomposable except if $C$ has a least element; in this latter case $\{(m, 0), (m,1)\}$ (here $m$ is the least element de $C$) is an acyclic component, and $K_{[C]}$ minus the vertex $(m,0)$ is isomorphic to $\check{K}_{[C]}$. 
\end{lemma}
 \begin{proof}
 Assertions (i). Every pair of distinct vertices of $C_3.C$ is included  in a $3$-cycle or a diamond.  Thus from Lemma \ref{separation}, $C_3.C$ is acyclically indecomposable. The second part of the sentence is obvious. Assertions (ii), (iii) and (iv)  follow directly from FactÊ\ref{fact:VTU} and the fact that $V_{2h+1}$, $T_{2h+1}$ and  $U_{2h+1}$ are indecomposable. Assertion (v)  follows by inspection. Note that every pair of distinct vertices of $H_{[C]}$ is included in a $3$-cycle. 
  Assertion (vi). The first part follows from the fact that $A'\times 2$ is an automous subset of $K_{[C]}$ for every initial interval $A'$ of $C$.  The second  part  follows from the fact that  every pair of distinct vertices of $K_{[C]}$ or of $K_{[C]}$ minus the vertex $(m,0)$ if $C$ has a least element $m$, is included in a $3$-cycle.    \end{proof}
  
From this, 
we deduce first:
    \begin{corollary}\label{acyclically indecomposable}
  Each member of $\mathfrak B$ except $K_{[\omega]}$ is acyclically indecomposable.  The tournament  $\check {K}_{[\omega]}$ is isomorphic to $K_{[\omega]}$ minus the vertex $(0,0)$.  In particular, it contains an isomorphic copy of $K_{[\omega]}$.
  \end{corollary}
  From our definitions, we have immediately this:
 
 \begin{fact}\label{incomp1}
The age of $V_{\omega}$  and   $V_{\omega^*}$ is  $\mathfrak V_{<\omega}$. The age of 
 $C_3.\omega$ and  $C_3.\omega^*$ is $\mathfrak C_{3,  <\omega}$.The age of    
$T_{\omega}$  and $T_{\omega^*}$ is $\mathfrak T_{<\omega}$. The age of  $U_{\omega}$ and  $U_{\omega^*}$ is $\mathfrak U_{<\omega}$. The age of  $H_{\omega}$ and  $H_{\omega^*}$ is $\mathfrak H_{<\omega}$. The age of  $K_{\omega}$ and  $K_{\omega^*}$ is $\mathfrak K_{<\omega}$.
\end{fact}

With the help of Lemma \ref {lem:twelveincec} we obtain:
 \begin{corollary} \label{increasingunion}For every $X_{[\alpha]}\in \mathfrak B$, $\check {X}_{[\alpha]}$ is an increasing union of $\check  {X}_{[\underline n]}$ for $n\in \N$. In particular, the age of $\check {X}_{[\alpha]}$ is the collection of finite tournaments which are embeddable in some $\check {X}_{[\underline n]}$ for some integer $n$.
\end{corollary}

\begin{lemma}\label{incomp2}
The six ages   $\mathfrak C_{3,  <\omega}$, $\mathfrak V_{<\omega}$,  $\mathfrak T_{<\omega}$, $\mathfrak U_{<\omega}$, $\mathfrak H_{<\omega}$, $\mathfrak K_{<\omega}$ are incomparable with respect to inclusion. 
\end{lemma}
\begin{proof}
 Let $\mathcal A:= \{\mathfrak C_{3,  <\omega}, \mathfrak V_{<\omega}, \mathfrak T_{<\omega}, \mathfrak U_{<\omega}, \mathfrak H_{<\omega}, \mathfrak K_{<\omega} \}$. Denote by $\neg \mathfrak C_{3,  <\omega}$ the set $\bigcup(\mathcal A\setminus \{\mathfrak C_{3,  <\omega}\})$ and define similarly $\neg \mathfrak V_{<\omega}$, $\neg \mathfrak T_{<\omega}$ etc. Let  $\tau_{1}$
(resp.$\tau_{2}$)be  a tournament obtained by replacing every vertex
of a 2-element tournament (resp. a
vertex of a 3-cycle) by a 3-cycle. 
We prove successively that 
 (i) $\tau_{1}\in  \mathfrak C_{3,  <\omega}\setminus \neg \mathfrak C_{3,  <\omega}$; (ii) $\tau_{2}\in  \mathfrak K_{<\omega}\setminus \neg \mathfrak K_{<\omega}$.
(iii) $T_{5}\in  \mathfrak T_{<\omega}\setminus \neg \mathfrak T_{<\omega}$; $V_{7}\in  \mathfrak V_{<\omega}\setminus \neg  \mathfrak V_{<\omega}$; (iv) $U_{7}\in  \mathfrak U_{ <\omega}\setminus \neg  \mathfrak U_{ <\omega}$; (v) $H_{[\underline 3]}\in \mathfrak H_{<\omega}\setminus\neg \mathfrak H_{<\omega}$;  

The proofs of the first and second  assertions  are immediate. Concerning the next three one, we may note that according to Lemma \ref{lem:TUV} and Corollary \ref{increasingunion}, $\mathfrak V_{<\omega}$,  $\mathfrak T_{<\omega}$, and $\mathfrak U_{<\omega}$ are pairwise incomparable w.r.t. inclusion. In fact, we derive these three assertions from the following observations:
\begin{enumerate}[{-}]
\item  The 3-cycle is the only  indecomposable subtournament of the
tournaments $C_{3}.\omega$ and $K_{[\omega]}$.
\item Up to
isomorphism,  $T_{2p+1}$
(resp. $V_{2p+1}$, $U_{2p+1}$) where $p \geq 2$, are the only
finite indecomposable subtournaments on at least 5 vertices   of $T_{[\omega]}$(resp. $V_{[\omega]}$, $U_{[\omega]}$). 
\item The tournaments $T_{5}$, $V_{7}$ and $U_{7}$ are not embeddable into $H_{[\omega]}$. 
\end{enumerate}
For the last assertion, we observe that the tournament 
$H_{[\underline 3]}$ is indecomposable and use the previous observations.
\end{proof}

\begin{lemma}\label{incomp}
Members of $\mathfrak B$ are pairwise incomparable with respect to embeddability.
\end{lemma}
\begin{proof}
According to Fact \ref{incomp1} and Lemma \ref{incomp2} it suffices to prove that: 
\begin{claim}\label{claim:incomp} If $\alpha \in \{\omega, \omega^*\}$, $
X_{[\alpha]}$ does not embed into $X_{[\alpha^*]}$.  
\end{claim}  If $X_{[\alpha]}$ is $C_{3[\alpha]}$, $V_{[\alpha]}$,   $T_{[\alpha]}$ or $H_{[\alpha]}$ this is obvious: $\alpha$ is embeddable in $X_{[\alpha]}$ but not in $X_{[\alpha^*]}$. If $X_{[\alpha]}= U_{[\alpha]}$, note that  for an arbitrary acyclic tournament $C$,  $U_{[C]}$ can be divided into two acyclic subsets $A_0$ and $A_1$ such that no $3$-cycle contains more than one vertex of $A_0$ and every pair of distinct vertices of $A_1$ is included in to some   $3$-cycle (set $A_0:=A\times \{0\}$ and $A_1:=A\times \{1\}$). Since in $U_{[\omega]}$, $A_0$ has type $\omega$, whereas in $U_{[\omega^*]}$, $A_0$ has type $\omega^*$, $U_{[\omega]}$ is not embeddable in $U_{[\omega^*]}$. If  $X_{[\alpha]}= K_{[\alpha]}$, note that
each autonomous set of  $K_{[\omega]}$ is finite whereas each autonomous set of  $K_{[\omega^*]}$ is cofinite. Hence, $K_{[\omega]}$ is not embeddable in $K_{[\omega^*]}$. 
\end{proof}

This is  the first part of Lemma  \ref{exponential}. We give the proof of the second part in the next section. 
\section{Profiles of members of $\mathfrak B$}\label {section:computprofile}
According to Fact \ref{incomp1}, our twelve acyclically indecomposable tournaments yield only six ages, those of $C_{3 {[\omega]}}$, $V_{[\omega]}$, $T_{[\omega]}$, $U_{[\omega]}$,  $H_{[\omega]}$ and  $K_{[\omega]}$.  For three of these tournaments, the exact values of the profile are know or easy to compute. For the others, we make no attempt  of an exact computation.

\subsection{Profile of $C_{3[\omega]}$}
The first values are 
$$1,1,1,2,3,4,6,9,13,19,28,41,60,88,129,189.$$ The sequence is  A000930 in \cite{sloane}. It
 satisfies the following recurrence $\varphi_{C_{3[\omega]}}(n) = \varphi_{C_{3[\omega]}}(n-1) + \varphi_{C_{3[\omega]}}(n-3)$ for $n\geq 3$. The Hibert series is 
 $H_{\varphi_{C_{3[\omega]}}}(x):=1/(1-x-x^3)$. According to \cite{cloitre},  $\varphi_{C_{3[\omega]}}(n) =\lfloor{ d*c^n + 1/2}\rfloor$ where 
 $c$  is the real root of $x^3-x^2-1$ and $d$ is the
              real root of $31*x^3-31*x^2+9*x-1$ ($c=1.465571231876768...$ and $d=
              0.611491991950812...$).

\subsection{Profile of $V_{[\omega]}$}
The first values are: $$1,1, 1, 2, 4, 9, 21, 48.$$
              
 \begin{fact} $\varphi_{V_{[\omega]}}(n)\geq 2^{n-5}$.
 \end{fact}
 \begin{proof} In $V_{[\omega]}$,  the vertex $a$ (or more exactly $\{a\}$) is the intersection of two $3$-cycles. We prove that the number of $n$-element restrictions of $V_{[\omega]}$ for which $a$ is the intersection of two $3$-cycles is at least  $2^{n-5}$. For that, let $n\geq 5$ be  an  integer.  For each subset $A$ of $\{1,\dots, n-5\}$, set $\overline A:=\{a\}\cup ((A\cup \{0, n-4\})\times \{0\}) \cup ((\{0,\dots n-4\}\setminus A)\times \{1\})$. As it is easy to see, the restrictions of $V_{[\omega]}$ to the $2^{n-5}$ subsets $A$ of  $\{1, \dots n-5\}$ are pairwise non isomorphic.  \end{proof}
              
\subsection{Profile of $T_{[\omega]}$}          
The tournament $T_{[\omega]}$ is diamond-free. Its age is the collection of finite diamond-free tournaments. There is a countable homogeneous tournament $L$ whose age is this collection. Thus $T_{[\omega]}$ and $L$ have the same profile. According to Cameron \cite{cameron2}:

\begin{equation}\label{eq:profile t}\varphi_{L}(n)= \frac{1}{2n}\sum_{d\vert n, d \text{odd}}\phi(d)2^{n/d}
\end{equation}

where $\phi$ is the Euler's totient function. 

As an immediate corollary we have $T_{[\omega]}(n)\geq {(2-\epsilon)}^{n}$.

\subsection{Profile of $U_{[\omega]}$}

\begin{lemma}\label{lem:profileu}$\varphi_{U_{[\omega]}}(n) \geq  (1-\epsilon)2^{n-2}$ for every  $\epsilon>0$ and  $n$ large enough. 
\end{lemma}
\begin{proof}
Let $h$ and $n$ be two integers with $5\leq 2h+1\leq n$.  Denote by $\Sigma_{h}(n)$ be the collection of 
 tournaments on  $n$ vertices which
 are a lexicographical sum of non-empty acyclic tournaments indexed by
 $U_{2h+1}$. Let $\Sigma (n):=\bigcup \{\Sigma_{h}(n): 5\leq 2h+1\leq n\}$ and   let $N_{h}(n)$, resp. $N(n)$, be the number of members of $\Sigma_{h}(n)$, resp. $\Sigma (n)$. 
 
  It is easy to check that each member of $\Sigma (n)$ is embeddable in   $U_{[\omega]}$. We claim that two members of $\Sigma (n)$ are isomorphic if and only if they are equal. Indeed, first
$\Sigma_{h}(n)$ and $\Sigma_{h'}(n)$ are disjoint whenever $h\not = h'$ (If they are not disjoint, then   since $U_{2h+1}$ and $U_{2h'+1}$ are indecomposable there are isomorphic, hence $h=h'$). Next, observe that  $U_{2h+1}$ is \emph{rigid}, that is the identity map is the unique automorphism of $U_{2h+1}$.  This observation follows readily from the fact that $\{h+1,\ldots,2h\}$ is the set of
centers of diamonds of $U_{2h+1}$ \cite{boudabbous}. From the rigidity of $U_{2h+1}$  follows that  a  lexicographic sum $\sum
 _{i\in U_{2h+1}} \underline {m}_i$ of non-empty acyclic tournaments $\underline {m}_i$ determines entirely the sequence $m_0, \dots, m_{2h}$.This proves our claim.  
 From this claim, we deduce first that:
$\varphi_{U_{[\omega]}}(n)\geq
 \sum\limits
 _{h=2}^{p} N_{h}(n)$ where $ p = \lfloor\frac{n-1}{2}\rfloor$.
We deduce next that  $N_{h}(n)$ is the number of integer 
 solutions of the equation: $n_{1}+\cdots+n_{2h+1} = n-2h-1$, that is $N_{h}(n)={ n-1\choose 2h}$. Combining these two facts, we have $\varphi_{U_{[\omega]}}(n) \geq \sum\limits
 _{h=2}^{p} { n-1\choose 2h} =  \sum\limits
_{h=0}^{p} {n-1\choose 2h}- 1 - {n-1\choose 2}=2^{n-2}- 1 - {n-1\choose 2}$. This proves the lemma.\end{proof}

\subsection{Profile of $H_{[\omega]}$}

\begin{lemma} \label{lem:profileh}$\varphi_{H_{[\omega]}}(n) \geq (1-\epsilon)2^{n-4}$ for every  $\epsilon >0$ and $n$ large enough. 
\end{lemma}
\begin{proof}
The proof is somewhat  similar to the proof of Lemma \ref{lem:profileu}.

Let $h\geq 3$ be an integer. Set $A_h:= \{(3k,i):k<h\}$ and $Z_h:=\{(3k+1,0), (3k+2, 1):k<h\}$. Let $Z\subseteq Z_h$; set $H_h(Z):= H_{[\omega]\restriction V_h\cup Z}$. 
Let $n$ be an integer with $n\geq 2h$.   Denote by $\Sigma'_{h}(Z, n)$  the collection of 
 tournaments $T$ on  $n$ vertices which
 are a lexicographical sum $\sum
 _{i\in H_{h}(Z)} \underline {m}_i$ of non-empty acyclic tournaments $\underline {m}_i$,  subject to the requirement that $m_i=1$ for each $i\not \in Z$. Let $\Sigma'_{h}(n):= \bigcup \{\Sigma'_{h}(Z, n): Z\subseteq Z_h\}$, $\Sigma'(n):= \bigcup \{\Sigma'_{h}( n): 6\leq 2h\leq n\}$ and  
 let $N'_{h}( Z, n)$,  resp. $N'_{h}( n)$, resp. $N'( n)$be the size of  $\Sigma'_{h}(Z, n)$, resp.  $\Sigma'_{h}(n)$, resp. $\Sigma'(n)$.   
 
  It is easy to check that each member of $\Sigma'(n)$ is embeddable in   $H_{[\omega]}$. We claim that two members of $\Sigma' (n)$ are isomorphic if and only if they are equal. This claim follows from the fact that  $H_h(Z)$ is indecomposable and rigid for every $Z \subseteq Z_h$. Indeed, note that from this fact $H_h(Z)$ and $H_h'(Z')$ are isomorphic if and only if there are equal, in particular $\Sigma'_{h}(n)$ and $\Sigma'_{h'}(n)$ are disjoint whenever $h\not = h'$. We leave the checking of the fact mentionned above to the reader   (we only note that $H_h(\emptyset)$ is isomorphic to $H_{\underline h}$.
From this claim, we deduce first that
$\varphi_{H_{[\omega]}}(n)\geq
 \sum\limits
 _{h=2}^{p} N_{h}(n)$ where $ p = \lfloor\frac{n}{2}\rfloor$.
We deduce next that  $N_{h}(n)$ is the number of integer 
 solutions of the equation: $n_{1}+\cdots+n_{2h-2} = n-2h$, that is $N_{h}(n)={n-3\choose 2h-3}$. Combining these two facts, we have $\varphi_{H_{[\omega]}}(n) \geq \sum\limits
 _{h=3}^{p} { n-3\choose 2h-3} =  \sum\limits
_{j=1}^{ \lfloor\frac{n-4}{2}\rfloor } {n-3\choose 2j+1}\geq 2^{n-4} - {n-3\choose 1}-1$. The conclusion of the lemma follows. \end{proof} 

\subsection{Profile of $K_{[\omega]}$.}
\begin{lemma}$\varphi_{K_{[\omega]}}(n) = 2^{n-2}$ for every  $n \geq 2$. 
\end{lemma}
\begin{proof}
We have $\varphi_{K_{[\omega]}}(0)=\varphi_{K_{[\omega]}}(1) =  \varphi_{K_{[\omega]}}(2)=1$.  We prove that:
\begin{equation}\label{eq:Komega}                     \varphi_{K_{[\omega]}}(n) = 1+ \sum\limits
 _{j=1}^{n-2} (n-j-1)\varphi_{K_{[\omega]}}(j)
 \end{equation}
 for $n\geq 3$. The lemma follows by induction on $n$. Since $\varphi_{K_{[\omega]}}(3)=2$,  formula (\ref{eq:Komega}) holds. Hence we may suppose $n\geq 4$.   
 
 Denote by $f_{K_{[\omega]}}(n)$ (resp.
$g_{K_{[\omega]}}(n)$) the number of strongly connected  (resp. non strongly connected) subtournaments of $K_{[\omega]}$ having $n$ vertices, these tournaments being counted up to
isomorphism. We have $f_{K_{[\omega]}}(0)=f_{K_{[\omega]}}(1)=f_{K_{[\omega]}}(2)=0$,  $f_{K_{[\omega]}}(3)=1$. More  generally, we have  $f_{K_{[\omega]}}(n)= \varphi_{K_{[\omega]}}(n-2)$ for $n\geq 4$. Indeed, every strongly connected subtournament of
$K_{[\omega]}$  having $n$ vertices, is obtained by dilating some vertex of
a $3$-cycle by a subtournament of $K_{[\omega]}$ having $(n-2)$ vertices of $K_{[\omega]}$. On an other hand, $g_{K_{[\omega]}}(n) = 1+ \sum\limits
 _{p=3}^{n-1} (n-p+1)f_{K_{[\omega]}}(p)$ for $n\geq 4$.
Indeed,  every non acyclic and non strongly connected subtournament  of $K_{[\omega]}$
has exactly one  strongly connected component which is not a singleton. Hence,  the number of non strongly
connected subtournaments of $K_{[\omega]}$ on $n$ vertices having a strongly connected component on $p$ vertices of a given isomorphy
type  is the number of integer solutions of the equation: $n_{1}+
n_{2}=n-p$. This number being  ${{n-p+1}\choose {1}}= n-p+1$ the above formula follows.
With the fact that  $\varphi_{K_{[\omega]}}(n)=
 f_{K_{[\omega}}(n)+
 g_{K_{[\omega]}}(n)$,  this yield formula (\ref{eq:Komega}).
 \end{proof}

 \section {Proof of Theorem \ref{morphology}}\label{sectiontheoremmorphology}
   
Let $T$ be an infinite acyclically indecomposable  tournament. We  prove that some member of $\mathfrak B$ is embeddable in $T$. The first step is:
\noindent \begin{claim}\label{claim:basic}  $V(T)$ contains an infinite subset $A$  such that:

 \begin{enumerate}
 \item Either every pair of distinct elements of $A$ is included into a $3$-cycle of $T$.
 \item Or  $V(T)$  contains no infinite subset  whose  pairs of distinct elements are included into a 3-cycle of $T$ and either:  
\begin{enumerate}
\item every pair of distinct elements of $A$ is included into a diamond of $T$ or
\item  every pair of elements of $A$  forms the end-vertices of  some double diamond of $T$ but $A$ contains no infinite subset  whose pairs of distinct elements are included into some diamond. 
\end{enumerate}
\end{enumerate}
\end{claim}
 
\noindent {\bf Proof of Claim \ref{claim:basic}.}
 Suppose that neither (1) nor (2-a) holds. 
 Let $f:\N \rightarrow V$ a one-to-one map. We define successively three subsets $X_1, X_2,  X_3$ made of  pairs $\{n,m\}$ of $[\N]^2$,  such that $n<m$,  depending wether $\{f(n),f(m)\}$ is contained into:
 \begin{enumerate}[{(a)}]
\item  some $3$-cycle of  $T$;
\item some diamond;
 \item   the endpoint of a self dual double diamond.   
\end{enumerate}
 
According to Lemma \ref{separation}, $[\N]^2 =X_1\cup X_2\cup X_3$. Hence, from Ramsey's Theorem, there is an infinite subset $Y$ of $\N$  and $i\in \{1,2,3\}$ such that $[Y]^2\subseteq X_i$.  Set $A:= \{f(n): n\in Y\}$. With our supposition, Case (a) and (b) are impossible. Thus  $A$ satisfies condition (2-b)  as claimed. \endproof. 
 
  Next, we prove that in case (1) some member of $\mathfrak B\setminus \{V_{[\omega]}, V_{[\omega^{*}]}\}$ is embeddable in $T$. In case (2-a), some member  of $\{C_{3 [\omega]}, C_{3[\omega^*]}, V_{[\omega]}, V_{[\omega^{*}]}\}$ is embeddable in $T$ and, in case (2-b), $C_{3 [\omega]}$ or  $C_{3[\omega^*]}$ is embeddable in $T$. The following lemmas take care of each case.

  \begin{lemma}\label{lem1}If a tournament $T$ contains an infinite subset $A$ such that every pair of distinct elements of $A$ is included into a $3$-cycle then  some member of $\mathfrak B\setminus \{V_{[\omega]}, V_{[\omega^{*}]}\}$ is embeddable in $T$.  
  \end{lemma}
  
  \begin{proof}
  Let $f:\N\rightarrow A$ be a one-to-one map. The hypothesis on $A$ allows to define  a map $g: [\N]^2\rightarrow V(T)$ such that $\{f(n),f(m), g(n, m)\}$ is a $3$-cycle of $T$ for every $n<m$. Let $\Phi:= \{f, g \}$,  let $\mathfrak L:=  <\omega, T, \Phi>$ and for a subset $X$ of $\N$, let $ \Phi_{\restriction X}:= \{f_{\restriction X}, g_{\restriction [X]^2}\}$ and let $\mathfrak L_{\restriction X}:= <\omega_{\restriction X}, T, \Phi_{\restriction X}>$. According to Theorem \ref{invar1}, there is an infinite subset $X$ of $\N$ such that $\mathfrak L_{\restriction X}$ is invariant. 
  
 Via a relabelling of $X$ with the integers, we may suppose that $X= \N$. Hence  $\mathfrak L$ is invariant. Let $A'$ be the image of $f$. 
 
\noindent \begin{claim}\label{claim:lem1}
\begin{enumerate}[{1.}]
\item $T(f(n), f(m))$ is constant on pairs $(n,m)$ such that $n<m$.
\item $T(f(n), g(m,k))$ is constant on triples $(n,m,k)$ such that $n<m<k$. 
\item   $g(n,m)\not \in \{f(k), g(n',m')\}$ for all $k, n<m<n'<m'$.
\item  If $T(g(n,m), g(m,n'))\not =T(g(n,m), g(n',m'))$ for some $n<m<n'<m'$ then for $D:= \{g(4k+i, 4k+i+1): k\in \N, i\in \{0, 1, 2\} \}$, $T_{\restriction D}$ is isomorphic to $C_3.\omega$ or to $C_3.\omega^*$.  
 \end{enumerate}
\end{claim} 
 \noindent {\bf Proof of Claim \ref{claim:lem1}.}
 Item 1. Since  $\mathfrak L$ is invariant, if $T(f(n), f(m))=1$ for some $n<m$ then $T(f(n'), f(m'))$=1 for all $n'<m'$. 
 
 \noindent Item 2. Same argument that in Item 1.
 
 \noindent Item 3. According to Item 1, if $T(f(n), f(m))=1$, $f$ is an isomorphism from $\omega$ onto $T_{\restriction A'}$. Similarly, if $T(f(n),f(m))=0$ for some $n<m$ then $f$ is an isomorphism from $\omega$ onto $T^*_{\restriction A'}$. In both cases, $T_{\restriction A'}$ is acyclic, hence it cannot contain a 
 $3$-cycle. This proves that there are no $k$ and  no $ n<m$ such that $g(n,m) = f(k)$.  Now,  suppose that  $g(n,m)=g(n',m')$ for some $n<m<n'<m'$; pick $m''$ with $m'<m''$. Let $h$ be the local isomorphism from $\omega$ to $\omega$ defined by $h(n)=n$, $h(m)=m$, $h(n')=m'$, $h(m')=m''$. Since $\mathfrak L$ is invariant, $g(n,m)=g(m', m'')$, hence  $g(n',m')=g(m', m'')$. But, as we have seen above, $T(f(n'), f(m'))=T(f(m'), f(m''))$. Since $\{f(n'), f(m'), g(n',m')\}$ is a $3$-cycle,  $T(f(n'), f(m'))=T(f(m'), g(n', m'))$. A similar argument yields  $T(g(m', m''), f(m'))=T(f(m'), f(m''))$, hence $T(g(m', m''), f(m'))=T(f(m'), g(n', m'))$. Since $T$ is a tournament, $g(n,m)\not= g(n',m')$, a contradiction. 
 
\noindent Item 4. Let $k\in \N$. Set $x_{i,k}:= g(4k+i, 4k+i+1)$ for $i\in \{0,1,2\}$ and $k\in \N$ and set $D_k:= \{x_{i,k}: i\in \{0,1,2\}\}$. Since $\mathfrak L$ is invariant, we have: $$T(x_{i,k}, x_{i+1, k})\not =T(x_{i,k}, x_{i+2, k})$$ 
  Again by the invariance of $\mathfrak L$, we have:  $$T(x_{i, k},  x_{i+1,k}) = T(  x_{i+1,k}, x_{i+2,k})$$ hence $T_{\restriction D_k}$ is a $3$-cycle.  
 By the invariance of $\mathfrak L$, $T(x_{i,k}, x_{i',k'})$ is constant on the pairs $(x_{i,k}, x_{i',k'})$ such that $k<k'$. If the value is $1$, $T_{\restriction D}$ is isomorphic to $C_3.\omega$ and if the value is $0$,  $T_{\restriction D}$ is isomorphic to $C_3.\omega^*$.  \endproof 
 
In order to get the conclusion of Lemma \ref{lem1}, we may suppose that neither  $C_3.\omega$, nor $C_3.\omega^*$, is embeddable in $T$. According to Item 1 of Claim \ref{claim:lem1}, $T(f(n),f(m))$ is constant on the pairs $(n,m)$ such that $n<m$. We may suppose that $T(f(n),f(m))=1$ (otherwise, it suffices to replace $T$ by $T^*$). 
We will consider two cases:

\noindent{\bf Case 1.} There are some $n<m<k$ such that $T( f(k), g(n, m))=1$. 

In this case, since $\mathfrak L$ is invariant, $T( f(k'), g(n',m'))=1$ for all $n'<m'<k'$. 

\noindent {\bf Case 2.}
$T( f(k),  g(n,m))=0$ for all  $n<m<k$.

Let $F:\N\times \{0, 1\}\rightarrow V(T)$ defined by setting $F(n,1):=f(n), F(n, 0):= g(n,n+1)$.
According to Item 3 of Claim \ref{claim:lem1}, $F$ is one to one. Let $T'$ be the tournament with vertex set $\N \times \{0, 1\}$, such that $T'(x, y)=T(F(x), F(y))$ for every pair of vertices of $\N \times \{0, 1\}$. 

\noindent \begin{claim}\label{claim2:lem1}
\begin{enumerate}[{1.}]
\item $T'((n, 1), (m,1))=1$ for $n<m
$. 
\item $T'((n, 0), (n+1, 0))=T'((n, 0), (m, 0))$ for $n<m$.
\item $T'((n,0), T'(n, 1))=1$
\item $T((n,1), (n+1,0))=T((n, 1), (m, 0))$ for $n<m$.
\end{enumerate}
\end{claim}
 \noindent {\bf Proof of Claim \ref{claim2:lem1}.}
\noindent Item1. $T'((n,1), T'(m, 1))=T(f(n), f(m))=1$.

\noindent Item 2. Since neither $C_{3[\omega]}$, nor $C_{3[\omega^*]}$, is embeddable in $T$,  we have $T'((n, 0), (n+1, 0))=T'((n, 0), (m, 0))$ for $n<m$.

\noindent Item 3. Since $\{f(n), f(n+1), g(n,n+1)\}$ is a $3$-cycle of $T$, $\{(n,1), (n+1, 1), (n, 0)\}$ is a $3$-cycle of $T'$. It follows that $T'((n,0), (n,1))=T'((n+1,1), (n,0))=1$. 

\noindent Item 4. Item 2 of Claim 1.
\endproof 

{\bf Suppose that Case 1 holds.} We have  $T'((m, 1), (n, 0))=1$ for all $m$, $n+1<m$.
Since $T'((n+1,1), (n,0))=1$,   we have $T'((n,1), (m,0))=1$ for all $n<m$. 
This added to Claim 2 insures that:
 
\noindent \begin{claim}\label{claim3:lem1}  $T'$ is a skew product of the $2$-element tournament $\underline 2$ by the chain $\omega$. \end{claim}
In order to conclude  the proof of Lemma \ref{lem1} there are four cases to consider.
Subcase 1.1. $T'((n, 1), (n+1, 0))=1$.
Subcase 1.1.1. $T'((n,0), (n+1,0))=1$.
In this case $T':= H_{[\omega]}$.
Subcase 1.1.2. $T'((n,0), (n+1, 0))=0$. In this case $T'$ is isomorphic to $K_{[\omega^*]}$. 
Subcase 1.2. $T'((n,1), (n+1,0))=0$. 
Subcase 1.2.1. $T( (n,0), (n+1, 0))=1$.
In this case $T'=T_{[\omega]}$.
Subcase 1.2.2. $T((n,0), (n+1, 0))=0$.
 In this case $T'$ is isomorphic to $U_{[\omega^*]}$.

 {\bf Suppose that Case 2 holds.} We have  $T'((m, 1), (n, 0))=0$ for all $m$, $n+1<m$.

\noindent \begin{claim}\label{claim4:lem1}   $T'((n,0), (n+1,0))=0$.
\end{claim}
 \noindent {\bf Proof of Claim \ref{claim4:lem1}.} Let $k\in \N$. Set $D_{k}:= \{(2k,0), (2k,1), (2k+1, 1))$. Set $D:= \cup\{D_k: k\in \N\}$. The tournament  $T'_{\restriction D_k}$ is a 3-cycle. If $T'((n,0), (n+1,0))=1$ then $T'((n,0), (m,0))=1$ for $n<m$. We have then $T'(x,y)=1$ whenever $x\in D_k, y\in D_k'$ , $k<k'$. Hence $T'_{\restriction D_k}$ is isomorphic to $C_3.\omega$. 
Contradicting our assumption. \endproof

\noindent \begin{claim}\label{claim5:lem1} Let $E':=\N\times 2\setminus \{(0,1)\}$. Then $T'_{\restriction E'}$ is isomorphic to $K_{[\omega]} $. \end{claim}

 \noindent {\bf Proof of Claim \ref{claim5:lem1}.} Let $G:\N\times 2\rightarrow \N\times 2$ defined by seting $G(n,0)= F(n, 1)$ and $G(n, 1)= F(n,0)$. Let $T''$ on $\N\times 2$ defined by $T''(x,y)=T'(x,y)$. One can easily check that $T''=K_{[\omega]}$. \endproof.
 
 With Claim \ref{claim5:lem1} the proof of Lemma \ref{lem1} is complete. 
\end{proof}

\begin{lemma}Let $T$ be  a tournament containing  no infinite subset whose  pairs of distinct elements are  included into a $3$-cycle. If  $T$ contains an infinite subset $A$ such that every  pair of distinct element of $A$ is included into a  diamond  then  
some member  of $\{C_3. \omega, C_3.\omega^*, V_{[\omega]}, V_{[\omega^{*}]}\}$ is embeddable in $T$.\end{lemma}
\begin{proof}
Let $f:\N\rightarrow A$ be a one-to-one map. We may define $f_i:[\N]^2\rightarrow V(T)$ for $i\in \{0,1\}$ so that for $n<m$, the set $\{f(n),f(m), f_i(n,m): i\in \{0,1\}\}$ forms a diamond and either $f(n)$ or $f(m)$ does not belong to the $3$-cycle of this diamond. 

 Let $\Phi:= \{f, f_i: i\in \{0,1\} \}$ and let $\mathfrak L:=  <\omega, T, \Phi>$.For a subset $X$ of $\N$, let $ \Phi_{\restriction X}:= \{f_{\restriction X}, f_{i \restriction [X]^2}: i\in \{0,1,\}\}$ and let $\mathfrak L_{\restriction X}:= <\omega_{\restriction X}, T, \Phi_{\restriction X}>$. According to Theorem \ref{invar1}, there is an infinite subset $X$ of $\N$ such that $\mathfrak L_{\restriction X}$ is invariant. 
 
Via a relabelling of $X$ with the integers, we may suppose that $X= \N$, that is   $\mathfrak L$   is invariant. 

 Since $\mathfrak L$ is invariant, $T(f(n),f(m))$ is constant on pairs $(n,m)$ such that $n<m$. With no loss of generality, we may suppose that 
 \begin{equation} \label{eq:f}T(f(n),f(m))=1. 
 \end{equation}

\noindent  {\bf Case 1}. Suppose that there is a pair $(n_0, m_0)$ such that $n_0<m_0$ and   $\{f(n_0), f_i(n_0,m_0): i\in \{0,1\}\}$ forms a $3$-cycle. With no loss of generality we may suppose that: $$T(f(n_0), f_1(n_0,m_0))= T(f_1(n_0,m_0),  f_0(n_0,m_0))= T(f_0(n_0, m_0), f(n_0))=1.$$  
 
\noindent \begin{claim}\label{claim1:lem2} 
Let $n<m<n'\leq m'$
\begin {enumerate}[{1.}]
\item $ T(f_i(m,n'), f(m'))=1$ for $i\in \{0, 1\}$. 
\item $T(f(n), f_1(m,n'))=1$.
\item $f_1(n,m)\not =f_1(n',m')$.
\item $f_i(n,m)\not =f(k)$ for $k\leq n$ or $m\leq k$. 
\end{enumerate}
\end{claim}
 \noindent {\bf Proof of Claim \ref{claim1:lem2}.}
\noindent Item 1. Since $\mathfrak L$ is invariant, $\{f(m), f_i(m,n'): i\in \{0,1\}\}$ is the $3$-cycle of $\{f(m') f(n'), f_i(m,n'): i\in \{0,1\}\}$. Hence $T(f_i(m,n'), f(n'))= T(f(m),f(n') )=1$. Now we may suppose  $m'>n'$. Since $T(f_i(m,n'), f(n'))=T(f(n'), f(m'))=1$ and $A$ does not contains a pair of distinct elements forming a $3$-cycle, $T(f_i(m,n'), f(m'))=1$. 

\noindent Item 2. We have  $T(f(n),f(m))=T(f(m),f_1(m,n'))=1$. Since $A$ does not contains a pair of distinct elements forming a $3$-cycle, 
$T(f(n),f_1(m,n'))=1$. 

\noindent Item 3. Suppose $f_1(n,m) =f_1(n',m')$. We have $T(f(n'), f_1(n',m'))=1$, whereas from Item 2, $ T(f_i(n,m), f(n'))=1$. Hence  $f_1(n,m) \not =f_1(n',m')$.

\noindent Item 4. If $k\geq m$, we have $T(f_i(n, m), f(k))$ hence the result. If $k=n$, $f_i(n,m)=f(n)$ is impossible by definition of $f_i$. If $k<n$ and  $f_i(n,m)=f(k)$ then select  $k'<k''<n''<m''$. By the invariance of $\mathfrak L$, get $f_i(n'',m'')=f(k')$ and $f_i(n'',m'')=f(k'')$, hence $f(k')=f(k'')$. A contradiction with the hypothesis that $f$ is one to one. \endproof

Note that from  Item 1 follows that one could choose $f_i(n,m)$ independent of $m$.

Let $F:\N\times \{0,1,2\}\rightarrow V(T)$ defined by $F(n,0)= f(3n)$, $F(n,i+1):=f_i(3n, 3n+1)$ for $i\in \{0,1\}$.

\noindent \begin{claim}\label{claim2:lem2} 
\begin{enumerate}[{1.}]
\item For $i,j\in \{0,1,2\}$, $T(F(n,i), (m,j))$ is constant on pairs $(n,m)$ such that $n<m$. 
\item $T(F(n,2), F(m,2))=1$ for $n<m$.
\item $T(F(n, 2), F(m, 1))=1 $ for $n<m$
\end{enumerate}
\end{claim}
 \noindent {\bf Proof of Claim \ref{claim2:lem2}.}
\noindent Item 1. $\mathfrak L$ is invariant. 

\noindent Item 2. If $T(F(n,2), F(m,2))=0$ then $T(f_1(3m, 3m+1), f_1(3n, 3n+1))=1$. Since $T(f(3m), f_1(3m,3m+1))=1$ and $T(f_1(3n,3n+1), f(3m))=1$, $\{ f_1(3m, 3m+1), f_1(3n, 3n+1), f(m)\}$ is a $3$-cycle. By  the invariance of $\mathfrak L$, every pair of distinct elements of $E':= \{f_1(3n', 3n'+1): n'\in \N\}$ contains a $3$-cycle, contradicting the  hypothese of the lemma.

\noindent Item 3. If $T(F(m,1), F(n, 2))=1$ then the set $\{F(m,1), F(n,2), F(m,2)\}$ forms a $3$-cycle and as in the item above we contradict the hypothesis in the lemma. 
\endproof. 

\noindent{\bf Subcase 1.1.}  $T(F(m,2), F(n, 1) )=1$ for some $n<m$.
Then this equality holds for all pairs $(n',m')$ such that $n'<m'$.
Let $a\not \in \N$ and  $G:\N\cup \{a\}\rightarrow V(T)$ defined by $G(a):= F(0,1)$, $G(n,0):=F(n, 0)$, $G(n,1):= F(n+1, 2)$. 

\noindent \begin{claim}\label{claim3:lem2} 
 $G$ is an embedding of $V_{[\omega]}$ into $T$. 
\end{claim}
\noindent {\bf Proof of Claim \ref{claim3:lem2}.}
Let $n\leq n'$. We have $T(G(n,0), G(n',0))=1$ from equation \ref{eq:f},  $T(G(n,1), G(n',1)=1$ from Item 2 of Claim\ref{claim2:lem2}, $T(G(n,0), G(n',1))=1$ from Item 2 and Item 3 of Claim\ref{claim1:lem2} and $T(G((n,1), G(n',0))=1$ from Item 1 of Claim\ref{claim1:lem2} if $n<n'$.  We have $T(G(a), G(n,0))=T(F(0,1), F(n,0))=T(F(n,1), F(n,0))=1$ from the hypothese of Case 1. And  we have $T(G(n,1), G(a))=T(F(n+1, 2),  F(0, 1))=1$. This proves our claim. 
\endproof

\noindent{\bf Subcase 1.2.}   $T(F(n,1), F(m, 2) )=1$ for every $n<m$.

\noindent \begin{claim}\label{claim4:lem2} 
 $F$ is one-to-one and  
$T(F(n, 1), F(m, 1)=1$ for every $n<m$.
\end{claim}
\noindent {\bf Proof of Claim \ref{claim4:lem2}.}
\noindent The first part of Claim \ref{claim4:lem2} is  obvious.If the second part does not hold, $\{F(n,1), F(m,2), F(m, 1)\}$ forms a $3$-cycle and every pair of distinct elements of $A':=\{F(n',1): n'\in \N\}$ is included into a $3$-cycle, contradicting the hypothesis of the lemma.
\endproof.

Let $G':\N\times 2 \cup \{a\}\rightarrow V(T)$ defined by $G'(a):= F(0,0)$, $G'(n,0):=F(n, 2)$, $G'(n, 1):= F(n+1, 1)$. 

\noindent \begin{claim}\label{claim5:lem2}  If   $T(F(m, 1), F(n, 0))=1$ for some $n<m$, $G'$ is an embedding from $ V_{[\omega]}$ into $T$. 
\end{claim}
Straightforward verification. 

Let $T'$ be defined on $\N\times \{0,1,2\}$ by $T'(x,y):=  T(F(x),F(y))$.
Let $T'_n:= T'_{\restriction \{n\}\times \N}$. 

\noindent \begin{claim}\label{claim5:lem2}  If $T(F(n, 0), F(m,1))=1$ for some $n<m$,  $T'$ is the $\omega$-sum of the $T_n$'s.
\end{claim}
Indeed, we have $T(F(n,i), F(m,j))=1$ for $n<m$.

\noindent {\bf Case 2}. Case 1 does not hold. In this case $\{f(m_0), f_i(n_0,m_0): i\in \{0,1\}\}$ forms a $3$-cycle. The treatment of this case is similar and left to the reader.  \end{proof}
 
 \begin{lemma} Let $T$ be  a tournament containing no infinite subset  whose pairs  of distinct vertices are  included into  a $3$-cycle. If $T$ contains an infinite subset $A$ such that every pair of distinct vertices of $A$ forms the end-vertices of some self-dual double diamond and is not included into a diamond  then either $C_{3[\omega]}$ or $C_{3[\omega^*]}$ is embeddable in $T$.
  \end{lemma} 
\begin{proof}
With Ramsey theorem, we may suppose that no pair of elements of $A$ is included into a $3$-cycle.  
Let $f:\N\rightarrow A$ be a one-to-one map. We may define $f_i:[\N]^2\rightarrow V(T)$ for $i\in \{0,1,2\}$ so that $f(n)$ and $f(m)$ are the end-vertices of the self-dual double diamond $\{f(n),f(m), f_i(n,m): i\in \{0,1,2\}\}$. 

 According to our construction, we have:
 
 \noindent \begin{claim}\label{claim1:lem3} \begin{enumerate}[{1.}]
 \item $T(f(n),f_i(n,m))=T(f(n),f(m))=T(f(_i(n,m),f(m))$.
 \item $T(f_0(n,m),f_1(n,m))=T(f_1(n,m),f_2(n,m))=T(f_2(n,m),f_0(n,m))$.
 \end{enumerate}
 \end{claim}
 Let $\Phi:= \{f , f_i: i\in \{0,1,2\} \}$ and let $\mathfrak L:=  <\omega, T, \Phi>$.For a subset $X$ of $\N$, let $ \Phi_{\restriction X}:= \{f_{\restriction X}, f_{i \restriction [X]^2}: i\in \{0,1,2\}\}$ and let $\mathfrak L_{\restriction X}:= <\omega_{\restriction X}, T, \Phi_{\restriction X}>$. According to Theorem \ref{invar1}, there is an infinite subset $X$ of $\N$ such that $\mathfrak L_{\restriction X}$ is invariant. 
  Via a relabelling of $X$ with the integers, we may suppose that $X= \N$, that is $\mathfrak L$ is invariant.

Let $F:\N\times \{0,1,2,3\}\rightarrow V(T)$ defined by $F(n,0)= f(3n)$, $F(n,i+1):=f_i(3n, 3n+1)$ for $i\in \{0,1,2\}$. 

\noindent \begin{claim}\label{claim2:lem3} $F$ is one-to-one. 
\end{claim}
\noindent {\bf Proof of Claim \ref{claim2:lem3}.}
Suppose first that $F_{i+1}(n)=F_{j+1}(m)$ for some $i,j\in\{0,1,2\}$, $n<m$.  
This means $f_i(3n,3n+1)=f_j(3m,3m+1)$. Since $\mathfrak L$ is invariant $f_i(3n,3n+1)=f_j(3m+1,3m+2)$. Hence $f_j(3m,3m+1)=f_j(3m+1,3m+2)$. From our construction, $T(f(3m), f(3m+1))=T(f_j(3m, 3m+1), f(3m+1))$ and $T(f(3m+1), f_j(3m+1, 3m+2))=T(f(3m+1), f(3m+2))$. Since $\mathfrak L$ is invariant, $T(f(n), f(m))$ is constant on pair $(n,m)$ such that $n<m$. In particular, $T(f(3m), f(3m+1))=T(f(3m+1), f(3m+2))$. Hence $T(f_j(3m,3m+1), f(3m+1))=T(f(3m+1), f_j(3m+1, 3m+2))$. Thus   $T(f_j(3m,3m+1), f(3m+1))=T(f(3m+1), f_j(3m, 3m+1))$, which is impossible since $f_j(3m,3m+1)$ and $ f(3m+1)$ are distinct.
Next, suppose that $F_{0}(n)=F_{i+1}(m)$ for some $n,m$, $n\not =m$. If  $n<m$, choose $n''<n'<m'$. Since $\mathfrak L$ is invariant, we have $f(3n'')= f_i(3m',3m'+1)$ and $f(3n')=f_i(3m',3m'+1)$ hence $f(3n'')=f(3n')$, contradicting the fact that $f$ is one-to-one. If $m<n$, choose $m'<n'<n''$ and use the same argument. \endproof

According to the above claim, we may define a tournament $T'$ 
with  vertex set $\N\times \{0,1,2,3\}$ such that $T'(x,y)= T(x,y)$. Indeed, 
let $T'_n:= T'_{ \restriction \{n\}\times \{0,1,2,3\}}$ for $n\in \N$:

\noindent \begin{claim}\label{claim3:lem3}$T'$ is the lexicographic  sum of the $T_n$'s, this sum being either an $\omega$-sum or an $\omega^*$-sum. 
\end{claim}
\noindent {\bf Proof of Claim \ref{claim3:lem3}.}
Let $i, j\in \{0,1,2,3\}$. Since $\mathfrak L$ is invariant, $T(F(n,i), F(m,j))$ is constant on pairs  $(n,m)$ such that $n<m$. Hence $T'$ is a skew product. We proceed directly. With no loss of generality, we may suppose $T(F(n,0), F(m,0))=1$ for $n<m$ (otherwise, replace $T$ by $T^*$), hence $T'((n,0), (m,0))=1$. According to Item 1 of Claim \ref{claim1:lem3}, we have $T(F(n,0), T(F(n, i+1))=1$. With this and the fact that  no pair of distinct elements of $A$ belong to a $3$-cycle, we also have $T(F(n,0), F(m, i+1))= T(F(n,i+1), F(m, 0))$ for all pairs $n,m$  such that $n<m$. Let  $\N_{i+1}:= \{(n, i+1): n\in \N\}$. If some pair of elements of $\N_{i+1}$ is included into  a $3$-cycle of $T'$, then all pairs are included into  a $3$-cycle. Hence, every pair of the infinite set $A_{i+1}:= \{ F(n, i+1): n\in \N\}$ would be included into some  $3$-cycle of elements of $T$, which is excluded. It follows that $T'((n,i+1), (m, i+1))=1$ if $n<m$. If $T'((n, i+1), (m,j+1))=0$ for some $n<m$, $i\not =j$, then the set $\{F(n,0),F(m,0), F(n, i+1), F(m,j+1)\}$ forms a diamond of $T'$ hence every pair of the infinite set $A_{0}:=\{F(n,0): n\in \N\}$ would be included into some diamond of $T$, which is excluded. From this $T'$ is the $\omega$-sum of the $T_n$'s, as claimed. 
\endproof

Since each $T_n$ contains a $3$-cycle, 
with Claim \ref{claim3:lem3} the proof of the lemma is complete.\end{proof}

\end{document}